\documentclass[11pt]{article}
\usepackage{amsmath,enumerate, amsfonts, amssymb, amsthm, graphicx, eepic, bbm}

\def\UseSection{
        \numberwithin{equation}{section}
    \theoremstyle{plain}
        \newtheorem{theorem}    {Theorem}[section]
        \DefineTheorems 
}

\def\DefineTheorems{
    
    \newtheorem{lemma}      [theorem] {Lemma}
    
    \newtheorem{prop}       [theorem] {Proposition}
    
    \newtheorem{cor}        [theorem] {Corollary}

    \theoremstyle{definition}
    \newtheorem{defn}       [theorem] {Definition}

    \theoremstyle{definition}

}

\newcommand{\bt}   {\begin{theorem}}
\newcommand{\et}   {\end  {theorem}}
\newcommand{\bl}   {\begin{lemma}}
\newcommand{\el}   {\end  {lemma}}
\newcommand{\bp}   {\begin{prop}}
\newcommand{\ep}   {\end  {prop}}
\newcommand{\bc}   {\begin{cor}}
\newcommand{\ec}   {\end  {cor}}
\newcommand{\bd}   {\begin{defn}}
\newcommand{\ed}   {\end  {defn}}

\newcommand{\ba}   {\begin{array}}
\newcommand{\ea}   {\end  {array}}
\newcommand{\be}   {\begin{enumerate}}
\newcommand{\ee}   {\end  {enumerate}}
\newcommand{\bi}   {\begin{itemize}}
\newcommand{\ei}   {\end  {itemize}}

\def\eq#1\en{\begin{equation}#1\end{equation}}
\def\eqsplit#1\ensplit{
    \begin{equation}\begin{split}#1\end{split}\end{equation}
    }
\def\eqalign#1\enalign{
    \begin{align}#1\end{align}
    }
\def\eqmul#1\enmul{
    \begin{multline}#1\end{multline}
    }
\newcommand{\eqarrstar} {\begin{eqnarray*}}
\newcommand{\enarrstar} {\end{eqnarray*}}
\newcommand{\eqarray}   {\begin{eqnarray}}
\newcommand{\enarray}   {\end{eqnarray}}
\newcommand{\nnb}   {\nonumber \\}

\newcommand{\lbeq}[1]  {\label{e:#1}}

\makeatletter
\newcommand{\labelcounter}[2]{{%
    \stepcounter{#1}
    \protected@write\@auxout{}%
    {\string\newlabel{#2}{{\csname the#1\endcsname}{\thepage}}}%
    {\ref{#2}}
    }}
\makeatother
\newcommand{\sss}   { \scriptscriptstyle }

\newcommand{\Rbold} {{\mathbb R}}

\newcommand{\Zbold} {{\mathbb Z}}

\newcommand{\Wcal}   {\mathcal{W}}

\newcommand{\chat}  {{ \hat{c}  }}

\newcommand{\Rd}    {{ {\Rbold}^d}}
\newcommand{\Zd}    {{ {\Zbold}^d }}

\newcommand{\spose}[1] {{\hbox to 0pt{#1\hss}} }
\newcommand{\ltapprox} {\mathrel{\spose{\lower 3pt\hbox{$\mathchar"218$}}
 \raise 2.0pt\hbox{$\mathchar"13C$}}}
\newcommand{\gtapprox} {\mathrel{\spose{\lower 3pt\hbox{$\mathchar"218$}}
 \raise 2.0pt\hbox{$\mathchar"13E$}}}

\newcommand{\ikx}  {{ i k \cdot x }}

\UseSection  
\newcommand{\R}     {\mathbb{R}}
\newcommand{\N}     {\mathbb{N}}

\newcommand\1{\mathbbm{1}}

\renewcommand{\d}{{\rm d}}

\newcommand{\eps}{\varepsilon}
\newcommand{\ssup}[1] {{\scriptscriptstyle{({#1})}}}

\newcommand{\Gk}{\hat{G}_z(k)}
\newcommand{\GN}{\hat{G}_z(0)}
\newcommand{\Dk}{\hat{D}(k)}
\newcommand{\dk}{\frac{\operatorname{d}\!k}{(2\pi)^d}}
\newcommand{\e}{\operatorname{e}}
\newcommand{\Td}{{\left[-\pi,\pi\right)^d}}

\newcommand{\ClK}{\hat{C}_{\lambda_z}(k)}
\newcommand{\ClN}{\hat{C}_{\lambda_z}(0)}

\newcommand{\G}{\hat{G}}

\newcommand{\const}{{\operatorname {const}\,}}
\newcommand{\Ul}{U_{\lambda_z}(k,l)}

\renewcommand{\phi}{\varphi}

\newcommand{\DL}{D}  

\newcommand{\am}{{\alpha\wedge2}}
\newcommand{\uv}{{\stackrel{\sss \rightharpoonup}{u}}}
\newcommand{\vv}{{\stackrel{\sss \rightharpoonup}{v}}} 
\newcommand{\bk}{{\boldsymbol{\rm k}}}
\newcommand{\bT}{{\boldsymbol{\rm T}}}
\newcommand{\bn}{{\boldsymbol{\rm n}}}
\newcommand{\bz}{{\boldsymbol{0}}}
\newcommand{\xir}{{\xi^{(r)}(n)}}

\newcommand{\cnT}[1]{{\hat{\boldsymbol{\rm c}}}^{\sss (#1)}}
\newcommand{\cnTl}[1]{\stackrel{\sss \le}{\boldsymbol{\rm c}}\!{}^{\sss (#1)}}
\newcommand{\cnTg}[1]{\stackrel{\sss >}{\boldsymbol{\rm c}}\!{}^{\sss (#1)}}

\begin{document}

\thispagestyle{empty}
\vspace{1cm}
\centerline{\LARGE \bf Long-range self-avoiding walk }
\vspace{0.2cm}
\centerline{\LARGE \bf converges to $\alpha$-stable processes}
\vspace{1cm}

\centerline {{\sc Markus Heydenreich}}
\vspace{.2cm}

\centerline{\em Vrije Universiteit Amsterdam, Department of Mathematics,}
\centerline{\em De Boelelaan 1081a, 1081 HV Amsterdam, The~Netherlands}
\centerline{\tt MO.Heydenreich@few.vu.nl}
\vspace{.4cm}

\centerline{\small(November 9, 2009)}
\vspace{.4cm}

\begin{quote}
  {\small {\bf Abstract:}}
  We consider a long-range version of self-avoiding walk in dimension $d>2(\am)$,
  where $d$ denotes dimension and $\alpha$ the power-law decay exponent of the coupling function.
  Under appropriate scaling we prove convergence to Brownian motion for $\alpha\ge2$, and to $\alpha$-stable L\'evy motion for $\alpha<2$.
  This complements results by Slade (1988), who proves convergence to Brownian motion for nearest-neighbor self-avoiding walk in high dimension.
\end{quote}

\vspace{1cm}


\vspace{1cm}
\hrule
\vspace{1em}

\section{Introduction and results}

\subsection{The model}
We study self-avoiding walk on the hypercubic lattice $\Zd$.
We consider $\Zd$ as a complete graph, i.e., the graph with vertex set $\Zd$ and corresponding edge set $\Zd\times\Zd$.
We assign each (undirected) bond $\{x,y\}$ a weight $D(x-y)$, where $D$ is a probability distribution specified in Section \ref{sectPropD} below. If $D(x-y)=0$, then we can omit the bond $\{x,y\}$.

\paragraph{Two-point function.}
For every lattice site $x\in\Zd$, we denote by
\begin{equation}\label{eqDefWn}
    \Wcal_n(x)=\{(w_0,\dots,w_{n}) \,|\, w_0=0,\, w_n=x,\, w_i\in\Zd, 1\le i\le n-1\}
\end{equation}
the set of $n$-step walks from the origin $0$ to $x$. We call such a walk $w\in\Wcal_n(x)$ \emph{self-avoiding} if $w_i\neq w_j$ for $i\neq j$ with $i,j\in\{0,\dots,n\}$. We define $c_0(x)=\delta_{0,x}$ and, for $n\ge1$,
\begin{equation}\label{eqDefCn}
    c_n(x):=\sum_{w\in\Wcal_n(x)}\prod_{i=1}^nD(w_i-w_{i-1})\,\1_{\text{\{$w$ is self-avoiding\}}}.
\end{equation}
where $D$ is specified below.
We refer to $D$ as the \emph{step} distribution, having in mind a random walker taking steps that are distributed according to $D$.
Without loss of generality we assume here that $D(0)=0$.

The self-avoiding walk measure is the measure $\mathbb Q_n$ on the set of $n$-step walks $\Wcal_n=\bigcup_{x\in\Zd}\Wcal_n(x)=\{0\}\times\mathbb{Z}^{dn}$
defined by
\begin{equation}\label{eqDefQn}
    \mathbb Q_n(w)
    :=\frac1{c_n}\,\prod_{i=1}^nD(w_i-w_{i-1})\,\1_{\text{\{$w$ is self-avoiding\}}},
\end{equation}
where $c_n=\sum_{x\in\Zd}c_n(x)$.

We consider the the Green's function $G_z(x)$, $x\in\Zd$, defined by
\begin{equation}\label{eqDefGsaw}
G_z(x)=\sum_{n=0}^\infty  c_n(x)\,z^n.
\end{equation}
We further introduce the \emph{susceptibility} as
\begin{equation}\label{eqDefChi}
\chi(z):=\sum_{x\in\Zd}G_z(x)
\end{equation}
and define $z_c$, the critical value of $z$, as the radius of convergence of the power series (\ref{eqDefChi}), i.e.\
\begin{equation}\label{eqDefZc}
    z_c:=\sup\,\left\{z\,|\,\chi(z)<\infty\right\}.
\end{equation}

The main part of our analysis is based on Fourier space analysis.
Unless specified otherwise, $k$ will always denote an arbitrary element from the Fourier dual of the discrete lattice, which is the torus $\Td$.
The Fourier transform of a function $f\colon\Zd\to\mathbb C$ is defined by
$    \hat f(k)=\sum_{x\in\Zd}f(x)\,\e^{ik\cdot x}$.

\paragraph{The step distribution $D$.}\label{sectPropD}
Let $h$ be a non-negative bounded function on $\Rd$ which is almost everywhere continuous, and symmetric under the lattice symmetries of reflection in coordinate hyperplanes and rotations by ninety degrees.
Furthermore we require $h$ to decay as $|x|^{-d-\alpha}$ as $|x|\to\infty$, where $\alpha>0$ is a parameter of the model.
In particular, there exists a positive constant $c_h$ such that
\begin{equation}\label{eqDefPowerLawDecay}
    h(x)\sim c_h\,|x|^{-d-\alpha}\qquad\text{whenever $|x|\to\infty$,}
\end{equation}
where $\sim$ denotes asymptotic equivalence, i.e., $f(x)\sim g(x)$ if $f(x)/g(x)\to1$.
For $\alpha\le2$ we make the stronger assumption that $h$ is completely rotation invariant on $\Rd$ (that is, not only by angles of 90 degrees as above).
Consequently, $\sum_{x\in\Zd} h(x/L)<\infty$ for all $L$, with $x/L=(x_1/L,\dots,x_d/L)$.

We then consider $D$ of the form
\begin{equation}\label{eqDefD}
    \DL(x)=\frac{h(x/L)}{\sum_{y\in\Zd}h(y/L)},
    \qquad x\in\Zd,
\end{equation}
where $L$ is a spread-out parameter (to be chosen large later on).
We note that the $\kappa$th moment $\sum_{x\in\Zd}|x|^\kappa \DL(x)$
does not exist if $\kappa\ge\alpha$, but exists and equals $O(L^\kappa)$ if $\kappa<\alpha$.

During the paper we shall make frequent use of the Landau symbols $O$ and $o$.
We denote $f=O(g)$ if $|f/g|$ is uniformly bounded. The bounding constant may depend on $d$, $\alpha$, $h$, but not on $n$, $k$, $z$, $u$, $\eps$ (these quantities are introduced later on). It may further depend on $L$ unless there is an explicit $L$-dependence in $g$ (like in the previous paragraph).
By $o(1)$ we denote terms that vanish as $n\to\infty$ (except for the appendix, where the limit $|k|\to0$ is considered).

\begin{lemma}[Properties of $D$]\label{lemmaD}
The step distribution $D$ satisfies the following properties:
\begin{enumerate}
    \item there is a constant $C$ such that, for all $L\ge 1$,
        \begin{equation}
        \label{eqPropD2}
        \|\DL\|_\infty\le CL^{-d};
        \end{equation}
    \item there is a constants $c>0$ such that
        \begin{eqnarray}
              \label{eqPropD3}
              1-\hat{D}(k)&>\hspace{0.2cm}c\qquad\qquad&\text{if $\|k\|_\infty\ge L^{-1}$,}\\
              \label{eqPropD4}
              1-\hat{D}(k)&<\hspace{0.2cm}2-c,\qquad &k\in\Td;
        \end{eqnarray}
    \item there is a constant $v_\alpha>0$ such that, as $|k|\to0$,
        \begin{equation}
              \label{eqPropD5}
              1-\hat{D}(k)\sim \;
              \begin{cases}
                v_\alpha |k|^\am\quad&\text{if $\alpha\neq2$},\\
                v_2|k|^2\log({1/\,|k|})\quad&\text{if $\alpha=2$}.
              \end{cases}
        \end{equation}
\end{enumerate}
\end{lemma}
\noindent
Chen and Sakai \cite[Prop.\ 1.1]{ChenSakai08} show that $D$ satisfies conditions (\ref{eqPropD2})--(\ref{eqPropD4}).
We prove in Appendix \ref{appendixD} that also (\ref{eqPropD5}) holds.
It follows from \cite[(1.7)]{ChenSakai08} that $v_\alpha\le O(L^{\am})$.

An example of $h$ satisfying all of the above is
\begin{equation}\label{eqDefHforPowerLaw}
    h(x)=(|x|\vee1)^{-d-\alpha},
\end{equation}
in which case $\DL$ has the form
\begin{equation}\label{eqExampleD}
\DL(x)=\frac{\left(|x/L|\vee1\right)^{-d-\alpha}}{\sum_{y\in\Zd}\left(|y/L|\vee1\right)^{-d-\alpha}},\qquad x\in\Zd.
\end{equation}


\subsection{Weak convergence of the end-to-end displacement.}
For $\alpha\in(0,\infty)$, we write
\begin{equation}\label{eqDefKn}
    k_n:=
    \begin{cases}
        k\,(v_\alpha n)^{-1/(\am)},\quad&\text{if $\alpha\neq2$}\\
        k\,(v_2 n\,\log\sqrt{n})^{-1/2},\quad&\text{if $\alpha=2$}
    \end{cases}
\end{equation}
so that
\begin{equation}\label{eqMain1}
    \lim_{n\to\infty}n\,[1-\hat D(k_n)]=|k|^\am.
\end{equation}

\bt[Weak convergence of end-to-end displacement]\label{thm-endpoint}
    Assume that $D$ is of the form \eqref{eqDefD},
    where the spread-out parameter $L$ is sufficiently large.
    Then self-avoiding walk in dimension $d>d_c=2(\am)$ satisfies
    \eq\lbeq{EndpointDistribution}
    \frac{\hat c_n(k_n)}{\hat c_n(0)}\to\exp\{-K_\alpha\,|k|^\am\}
    \qquad\text{as $n\to\infty$,}
    \en
    where
    \begin{equation}\label{eqDefXi}\displaystyle
        K_\alpha=
        \Big(1+\sum_{x\in\Zd}\sum_{n=2}^\infty n\,\pi_n(x)\,{z_c}^{n-1}\Big)^{-1}\times
        \begin{cases}
            1,\;&\text{if $\alpha\le2$};\\
            \displaystyle
            1+(2dz_cv_\alpha)^{-1}\sum_{x\in\Zd}\sum_{n=2}^\infty |x|^2\,\pi_n(x)\,{z_c}^n,\;&\text{if $\alpha>2$.}
        \end{cases}
    \end{equation}
\et
The quantities $\pi_n(x)$ appearing in \eqref{eqDefXi} are known as lace expansion coefficients.
We do not perform the lace expansion in this paper.
References to the derivation of the lace expansion and various bounds on these lace expansion coefficients are given later on.
Under the conditions of Theorem \ref{thm-endpoint},
\eqref{eqLemmaMoments} and \eqref{eqPiBd4} below imply that both sums appearing in \eqref{eqDefXi} are finite.
However, the quantities $\pi_n(x)$ are given in terms of an alternating sum, cf.\ \eqref{eqPiBd1}, and their sign is not known.
Nevertheless, both sums appearing in \eqref{eqDefXi} can be made smaller than 1 by taking $L$ large enough, as proven in \cite{HofstSlade02} for $\alpha>2$, and for $\alpha\le2$ it follows the lines of \cite[Section 6.2.2]{MadraSlade93} in combination with \cite{HeydeHofstSakai08}.
Consequently, $K_\alpha\in(0,\infty)$.

\subsection{Mean-$r$ displacement.}
The \emph{mean-$r$ displacement} is defined as
\begin{equation}\label{eqDefXiR}
    \xir:=\left(\frac{\sum_{x\in\Zd}|x|^rc_n(x)}{c_n}\right)^{1/r},
\end{equation}
where we recall $c_n=\sum_{x\in\Zd}c_n(x)=\hat c_n(0)$.
For $r=2$ this is the mean-square displacement, and already well understood.
For example, van der Hofstad and Slade \cite{HofstSlade02} prove the following rather general version:
\bt[Mean-square displacement {\cite[Theorem 1.1.b]{HofstSlade02}}]\label{thm-xi2}
Consider self-avoiding walk with step distribution $D$ given in Section \ref{sectPropD} with $\alpha>2$.
Then there exist constants $C>0$ and $\delta>0$ (both depending on $d,\alpha,h,L$) such that, as $n\to\infty$,
\begin{equation}\label{eqXi2}
    \frac{1}{c_n}\sum_{x\in\Zd}|x|^2c_n(x)
    =C\,n\,(1+O(n^{-\delta})).
\end{equation}
\et
The proof of Theorem \ref{thm-xi2} is also based on lace expansion.
In the sequel we prove a complementary result for $r<2$.
To this end, we write $f\asymp g$ if there are uniform positive constants with $cg\le f \le Cg$.
\bt[Mean-$r$ displacement]\label{thm-xiR}
Under the assumptions of Theorem \ref{thm-endpoint},
for any $r<\am$,
\begin{equation}\label{eqXiR}
    \xir\asymp
        \begin{cases}
            n^{1/(\am)},\quad&\text{if $\alpha\neq2$},\\
            (n\,\log n)^{1/2}
            ,\quad&\text{if $\alpha=2$},
        \end{cases}
\end{equation}
as $n\to\infty$.
\et
Recently, Chen and Sakai \cite{ChenSakai10} found the proof that \eqref{eqXiR} holds for \textit{all} $r\in(0,\alpha)$,
for long-range self-avoiding walk and long-range oriented percolation.

\subsection{Convergence to Brownian motion and $\alpha$-stable processes.}
In order to deal with the cases $\alpha=2$ and $\alpha\neq2$ simultaneously, we write
\begin{equation}\label{eqDefFAlpha}
    f_\alpha(n)
    =   \begin{cases}
            (v_\alpha n)^{-1/(\am)}\quad&\text{if $\alpha\neq2$},\\
            (v_2 n\,\log\sqrt{n})^{-1/2}\quad&\text{if $\alpha=2$},
        \end{cases}
\end{equation}
such that, for example, $k_n=f_\alpha(n)\,k$, cf.\ \eqref{eqDefKn}.
Given an $n$-step self-avoiding walk $w$, define
\begin{equation}\label{eqDefXn}
    X_n(t)=\left(2dK_\alpha\right)^{-\frac1\am}\,f_\alpha(n)\,w(\lfloor nt \rfloor),
    \qquad
    t\in[0,1].
\end{equation}
We aim to identify the scaling limit of $X_n$, and the appropriate space to study the limit is the space of $\R^d$-valued c\`adl\`ag-functions $D([0,1],\Rd)$ equipped with the Skorokhod topology.

For $\alpha\in(0,2]$, $W^{(\alpha)}$ denotes the standard $\alpha$-stable L\'evy measure, normalized such that
\begin{equation}
    \int \e^{ik\cdot B^{(\alpha)}(t)}\,\d W^{(\alpha)}=\e^{-|k|^\alpha t/(2d)},
\end{equation}
where $B^{(\alpha)}$ is a (c\`adl\`ag version of) standard symmetric $\alpha$-stable L\'evy motion (in the sense of \cite[Definition 3.1.3]{SamorTaqqu94}).
Note that $W^{(2)}$ is the Wiener measure, and $B^{(2)}$ is Brownian motion.
By $\langle\cdot\rangle_n$ we denote expectation with respect to the self-avoiding walk measure $\mathbb Q_n$ in \eqref{eqDefQn}.

\bt[Weak convergence to $\alpha$-stable processes and Brownian motion]\label{thm-pathconvergence}
    Under the assumptions in Theorem \ref{thm-endpoint},
    \begin{equation}
        \lim_{n\to\infty}\langle f(X_n)\rangle_n=\int f\,\d W^{(\am)},
    \end{equation}
    for every bounded continuous function $f\colon D([0,1],\Rd)\to\mathbb R$. That is to say,
    $X_n$ converges in distribution to an $\alpha$-stable L\'evy motion for $\alpha<2$,
    and to Brownian motion for $\alpha\ge2$.
    Equivalently, $\mathbb Q_n$ converges weakly to $W^{(\am)}$.
\et

In order to prove convergence in distribution as a process, we need two properties: $(i)$ the convergence of finite-dimensional distributions, and $(ii)$ tightness of the family $\{X_n\}$.
We shall now consider the former.

Convergence of finite-dimensional distributions means for every $N=1,2,3,\dots$, any $0<t_1<\cdots<t_N\le1$, and any bounded continuous function $g\colon\mathbb R^{dN}\to\mathbb R$,
\begin{equation}\label{eqFinDimConv2}
    \lim_{n\to\infty}\left\langle g\big(X_n(t_1),\dots,X_n(t_N)\big)\right\rangle_n
    =\int g\big(B^{(\am)}(t_1),\dots,B^{(\am)}(t_N)\big)\,\d W^{(\am)}.
\end{equation}
Convergence of characteristic functions determines convergence in distribution, it is therefore sufficient to consider functions $g$ of the form
\begin{equation}
    g(x_1,\dots,x_N)=\exp\{i\,\bk\cdot(x_1,\dots,x_N)\},
\end{equation}
where $\bk=\big(k^{\sss (1)},\dots,k^{\sss (N)}\big)\in\R^{dN}$ and $x_i\in\Rd$, $i=1,\dots,N$.
We rather use the equivalent form
\begin{equation}
    g(x_1,\dots,x_N)=\exp\{i\,\bk\cdot(x_1,x_2-x_1,\dots,x_N-x_{N-1})\},
\end{equation}
which better fits in our setting.

%
For $\bn=(n^{\sss (1)},\dots,n^{\sss (N)})\in\N^N$, with $n^{\sss (1)}<\dots<n^{\sss (N)}$,
we define
\begin{equation}\label{eqDefCnN}
    \begin{split}
    \cnT{N}_{\bn}(\bk):=\sum_{x_1,x_2,\dots,x_{n^{\sss (N)}}}&
      \exp\left\{i\sum_{j=1}^Nk^{\sss (j)}\cdot\left(x_{n^{\sss (j)}}-x_{n^{\sss (j-1)}}\right)\right\}\\
      &\times\prod_{i=1}^{n^{\sss (N)}}D(x_i-x_{i-1})\,
      \1_{\{(0,x_1,x_2,\dots,x_{n^{\sss (N)}}) \text{ is self-avoiding}\}}
    \end{split}
\end{equation}
as the $N$-dimensional version of the Fourier transform of \eqref{eqDefCn}, with $n^{\sss (0)}=0$.
An alternative representation is
\begin{equation}\label{eqDefCnN2}
    \cnT{N}_{\bn}(\bk)
      =\sum_{w\in \Wcal_{n^{\sss (N)}}}
      \e^{i\bk\cdot\Delta w(\bn)}
      W(w)\;
      \1_{\{\text{$w$ is self-avoiding}\}},
\end{equation}
where $W(w)=\prod_{i=1}^{|w|}D(w_i-w_{i-1})$ is the \emph{weight} of the walk $w$ ($|w|$ denotes the length) and
$$  \bk\cdot\Delta w(\bn)
    =\sum_{j=1}^Nk^{\sss (j)}\cdot\left(w_{n^{\sss (j)}}-w_{n^{\sss (j-1)}}\right).
 $$

We fix a sequence $b_n$ diverging to infinity slowly enough such that
\begin{equation}\label{eqDefBn}
    f_\alpha(n)^{\alpha\wedge1}\,b_n=o(1),
\end{equation}
for example $b_n=\log n$.
\bt[Finite-dimensional distributions]\label{prop-FinDimConv}
Let $N$ be a positive integer,
$k^{\sss (1)},\dots,$ $k^{\sss (N)}\in\R^d$,
$0=t^{\sss (0)}<t^{\sss (1)}<\dots<t^{\sss (N)}\in[0,1]$,
and $g=(g_n)$ a sequence of real numbers satisfying $0\le g_n\le b_n/n$.
Denote
$$  \bk_n
    =\big(k_n^{\sss (1)},\dots,k_n^{\sss (N)}\big)
    =f_\alpha(n)\,\big(k^{\sss (1)},\dots,k^{\sss (N)}\big),
$$
$$n\bT=\big(\lfloor nt^{\sss (1)}\rfloor,\dots,\lfloor nt^{\sss (N-1)}\rfloor,\lfloor nT\rfloor\big)$$
with $T=t^{\sss (N)}(1-g_n)$.
Under the conditions of Theorem \ref{thm-endpoint},
\begin{equation}\label{eqFinDimConv}
    \lim_{n\to\infty}\frac{\cnT{N}_{n\bT}(\bk_n)}{\chat_{nT}(0)}
    =\exp\left\{-K_\alpha\,\sum_{j=1}^N|k^{\sss (j)}|^\am\;(t^{\sss (j)}-t^{\sss (j-1)})\right\}
\end{equation}
holds uniformly in $g$.
\et

The presence of the sequence $g_n$ might appear unclear at this point,
it is there for a technical reason:
The proof of Theorem \ref{prop-FinDimConv} is carried out by induction over $N$ and some flexibility is needed in the endpoint.

Let us emphasize that \eqref{eqFinDimConv} has indeed the required form.
Let $k^{\sss (1)},\dots,$ $k^{\sss (N)}\in\R^d$ and
$0=t^{\sss (0)}<t^{\sss (1)}<\dots<t^{\sss (N)}\in[0,1]$ be given.
We apply Theorem \eqref{prop-FinDimConv} with $N+1$ and $g_n\equiv0$,
where $k^{\sss (N+1)}=0$ and $T=t^{\sss (N+1)}=1$, so that
$n\bT=\big(\lfloor nt^{\sss (1)}\rfloor,\dots,\lfloor nt^{\sss (N)}\rfloor,n\big)$. Then
\begin{eqnarray}
    \Big\langle\exp\big\{i\,\bk\cdot\Delta X_n(n\bT)
        \big\}\Big\rangle_{n}
    &=& \left\langle\exp\left\{i\left(2dK_\alpha\right)^{-\frac1\am}\bk_n
        \cdot\Delta \omega(n\bT)\right\}\right\rangle_{n}\nnb
    &=& \frac{\cnT{N}_{n\bT}\!\left(\left(2dK_\alpha\right)^{-\frac1\am}\bk_n\right)}{\chat_{n}(0)},
    \nonumber
\end{eqnarray}
and this converges to
$$ \exp\left\{-\frac1{2d}\,\sum_{j=1}^N|k^{\sss (j)}|^\am\;(t^{\sss (j)}-t^{\sss (j-1)})\right\}$$
as $n\to\infty$, as we aim to show for \eqref{eqFinDimConv2}.
Thus the finite dimensional distributions of (long-range) self-avoiding walk converge to those of an $\alpha$-stable L\'evy motion, which proves that this is the only possible scaling limit.

\subsection{Discussion and related work}
Long-range self-avoiding walk has rarely been studied.
Klein and Yang \cite{YangKlein88} show that the endpoint of a weakly self-avoiding walk jumping $m$ lattice sites \emph{along the coordinate axes} with probability proportional to $1/m^2$, is Cauchy distributed. 
A similar result for strictly self-avoiding walk is obtained by Cheng \cite{Cheng00}.

In a previous paper \cite{HeydeHofstSakai08} it is shown that long-range self-avoiding walk exhibits mean-field behavior above dimension $d_c=2(\am)$.
More specifically, it is shown that under the conditions of Theorem \ref{thm-endpoint}, the Fourier transform of the critical two-point function satisfies
$\hat G_{z_c}(k)=(1+O(\beta))/(1-\Dk)$, where $\beta=O(L^{-d})$ is an arbitrarily small quantity.
Hence, on the level of Fourier transforms, the critical two-point functions of long-range self-avoiding walk and long-range \emph{simple} random walk are very close.
Indeed, the results in \cite{HeydeHofstSakai08} suggest that the two models behave similarly for $d>d_c$, and we confirm this in a rather strong form by showing that both objects have the same scaling limit.

Chen and Sakai \cite{ChenSakai09} prove an analogue of Theorem \ref{thm-endpoint} for oriented percolation, and in fact our method of proving Theorem \ref{thm-endpoint} is very much inspired by the method in \cite{ChenSakai09}.
The bounds on the diagrams are different for the two different models, but the general strategy works equally well with either model.
In particular, the \emph{spatial} fractional derivatives as in \eqref{eqMom1} are used for the first time in \cite{ChenSakai09}.

Slade \cite{Slade88, Slade89} proves convergence of the \emph{nearest-neighbor} self-avoiding walk to Brownian motion in sufficiently high dimension, using a finite-memory cut-off.
Hara and Slade \cite{HaraSlade92a} provide an alternative argument by using fractional derivative estimates. An account of the latter approach is contained in the monograph
\cite[Sect.\ 6.6]{MadraSlade93}.
All of these proofs use the lace expansion, which was introduced by Brydges and Spencer \cite{BrydgSpenc85} to study weakly self-avoiding walk.

\section{The scaling limit of the endpoint: Proof of Theorem \ref{thm-endpoint}}
\subsection{Overview of proof}\label{sectLaceExpansion}\noindent
The lace expansion obtains an expansion of the form
\begin{equation}\label{eqExpansionX}
    {c}_{n+1}(x) = ({D}\ast{c}_n)(x)
    + \sum_{m=2}^{n+1} \left({\pi}_m \ast {c}_{n+1-m}\right)(x)
\end{equation}
for suitable coefficients $\pi_m(x)$,
see e.g.\ \cite[Sect.\ 2.2.1]{Hofst05} or \cite[Sect.\ 3]{Slade06} for a derivation of the lace expansion.
We multiply (\ref{eqExpansionX}) by $z^{n+1}$ and sum over $n\ge0$.
By letting
\begin{equation}\label{eqDefPi}
\Pi_z(x)=\sum_{m=2}^\infty\pi_m(x)z^m
\end{equation}
for $z\le z_c$, and recalling
$G_z(x)=\sum_{n=0}^\infty c_n(x)z^n$,
this yields
\begin{equation}\label{eqExpansionX2}
    G_z(x)=\delta_{0,x}+z(D\ast G_z)(x)+(G_z\ast\Pi_z)(x).
\end{equation}

We proceed by proving Theorem \ref{thm-endpoint} subject to certain bounds on the lace expansion coefficients $\pi_n(x)$ to be formulated below.
A Fourier transformation of (\ref{eqExpansionX2}) yields
\begin{equation}\label{eqExpansionK1}
    \Gk=1+z\,\Dk\, \Gk+\Gk\,\hat\Pi_z(k),\qquad k\in\Td,
\end{equation}
and this can be solved for $\Gk$ as
\begin{equation}\label{eqProof1}
    \Gk^{-1}=1-z\,\Dk-\hat\Pi_z(k),\qquad k\in\Td.
\end{equation}
Since $z_c$ is characterized by $\hat G_{z_c}(0)^{-1}=0$, one has $\hat \Pi_{z_c}(0)=1-z_c$, and hence
\begin{equation}\label{eqProof2}
    \Gk^{-1}=(z_c-z)\,\Dk
             +\left(\hat\Pi_{z_c}(k)-\hat\Pi_z(k)\right)
             +z_c(1-\Dk)
             +\left(\hat\Pi_{z_c}(0)-\hat\Pi_{z_c}(k)\right).
\end{equation}
If we let
\begin{eqnarray}
  \label{eqProof3}
  A(k) &:=&  \Dk+\partial_z\hat\Pi_z(k)\big|_{z=z_c},\\
  \label{eqProof4}
  B(k) &:=&  1-\Dk+\frac1{z_c}\left(\hat\Pi_{z_c}(0)-\hat\Pi_{z_c}(k)\right),\\
  \label{eqProof5}
  E_z(k) &:=& \frac{\hat\Pi_{z_c}(k)-\hat\Pi_{z}(k)}{z_c-z}-\partial_z\hat\Pi_z(k)\big|_{z=z_c},
\end{eqnarray}
then
\begin{eqnarray}
  \label{eqProof6}\nonumber
  z_c\,\Gk
  &=& \frac{1}{[1-z/z_c]\left(A(k)+E_z(k)\right)+B(k)} \\
  \label{eqProof7}
  &=& \frac{1}{[1-z/z_c]\,A(k)+B(k)} -\Theta_z(k),
\end{eqnarray}
where
\begin{equation}\label{eqProof8}
    \Theta_z(k)=\frac{[1-z/z_c]\,E_z(k)}
                    {\big([1-z/z_c]\left(A(k)+E_z(k)\right)+B(k)\big)\,
                     \big([1-z/z_c]\,A(k)+B(k)\big)}.
\end{equation}
If $\Gk^{-1}$ is understood as a function of $z$, then $A(k)$ denotes the linear contribution, $E_z(k)$ denotes the higher order contribution (which will turn out to be asymptotically negligible), and $B(k)$ denotes the constant term.
The denominators in \eqref{eqProof7}--\eqref{eqProof8} are positive for $z<z_c$, cf.\ \eqref{eqErrorProof14}--\eqref{eqErrorProof15} below.

For the first term in (\ref{eqProof7}) we write
\begin{equation}\label{eqProof9}
    \frac{1}{[1-z/z_c]\,A(k)+B(k)}
    =\frac{1}{A(k)+B(k)}\sum_{n=0}^\infty\left(\frac
     z{z_c}\right)^n\left(\frac{A(k)}{A(k)+B(k)}\right)^n,
\end{equation}
and the geometric sum converges whenever $z<z_c\left(A(k)+B(k)\right)/A(k)$; the latter term approximates $z_c$ as $|k|\to0$.
For $z<z_c$, we can write $\Theta_z(k)$ as a power series,
\begin{equation}\label{eqProof10}
    \Theta_z(k)=\sum_{n=0}^\infty\theta_n(k)\,z^n.
\end{equation}
Since $\Gk=\sum_{n=0}^\infty \hat c_n(k)z^n$ and $B(0)=0$, we thus obtained
\begin{equation}\label{eqProof11}
    \hat c_n(k)=\frac1{z_c}\left(\frac{z_c^{-n}}{A(k)+B(k)}\left(\frac{A(k)}{A(k)+B(k)}\right)^n-\theta_n(k)\right),
    \quad
    \hat c_n(0)=\frac1{z_c}\left(\frac{z_c^{-n}}{A(0)}-\theta_n(0)\right).
\end{equation}
In Section \ref{sectErrorBound} we prove the following bound on the error term $\theta_n$:
\bl\label{lemmaError}
    Under the conditions of Theorem \ref{thm-endpoint}, $|\theta_n(k)|\le O(z_c^{-n}\,n^{-\eps})$ for all $\eps\in\left(0,\big(\frac{d}{\am}-2\big)\wedge1\right)$
    uniformly in $k\in\Td$.
\el
Equation \eqref{eqProof11} and Lemma \ref{lemmaError} imply the following corollary:
\bc\label{corol-CnAsymptotic}
Under the conditions of Theorem \ref{thm-endpoint},
\begin{equation}\label{}
    \hat c_n(0)={\Xi}\,{z_c^{-n}}\left(1+O(n^{-\eps})\right),
\end{equation}
where $\eps\in\left(0,\big({d}/({\am})-2\big)\wedge1\right)$ and
\begin{equation}\label{}
    \Xi= \left[{z_c\,A(0)}\right]^{-1}
     = \left[z_c+\sum_{x\in\Zd}\sum_{m=2}^\infty m\,\pi_m(x)\, z_c^m\right]^{-1}
     \in(0,\infty).
\end{equation}
\ec

By (\ref{eqProof11}) and Lemma \ref{lemmaError},
for $\eps\in\left(0,\big(\frac{d}{\am}-2\big)\wedge1\right)$ an all $k\in\Rd$ such that $k_n\in[-\pi,\pi)^d$,
\begin{eqnarray}
    \frac{\hat c_n(k_n)}{\hat c_n(0)}
    &=&\big(1+O(n^{-\eps})\big)\frac{A(0)}{A(k_n)+B(k_n)}\left(\frac{A(k_n)}{A(k_n)+B(k_n)}\right)^n+O(n^{-\eps})\nnb
    \label{eqProof12}
    &=& \big(1+O(n^{-\eps})\big)\frac{A(0)}{A(k_n)+B(k_n)}\\
    &&  \quad\times\left(1+\frac{-n(1-\hat D(k_n))\,\big(A(k_n)+B(k_n)\big)^{-1}\,{B(k_n)}(1-\hat D(k_n))^{-1}}{n}\right)^n
        +O(n^{-\eps}).\nonumber
\end{eqnarray}
As $n\to\infty$, we have that
$n(1-\hat D(k_n))\to|k|^\am$ by (\ref{eqMain1}),
$$A(k_n)\to A(0)=1+\sum_{x\in\Zd}\sum_{m=2}^\infty m\,\pi_m(x)\, z_c^{m-1}.$$
The convergence
\begin{equation}\label{eqProof13}
    \lim_{n\to \infty}\frac{B(k_n)}{1-\hat D(k_n)}=
        \begin{cases}
            1,\quad&\text{if $\alpha\le2$};\\
            1+(2dz_cv_\alpha)^{-1}\sum_{x\in\Zd}|x|^2\,\Pi_{z_c}(x),\quad&\text{if $\alpha>2$.}
        \end{cases}
\end{equation}
follows directly from the following proposition:
\bp\label{propMain2}
Under the conditions of Theorem \ref{thm-endpoint},
\begin{equation}\label{eqMain2}
    \lim_{|k|\to 0}\frac{\hat\Pi_{z_c}(0)-\hat\Pi_{z_c}(k)}{1-\Dk}=
        \begin{cases}
            0,\quad&\text{if $\alpha\le2$};\\
            (2dv_\alpha)^{-1}\sum_{x\in\Zd}|x|^2\,\Pi_{z_c}(x),\quad&\text{if $\alpha>2$.}
        \end{cases}
\end{equation}
\ep
If a sequence $h_n$ converges to a limit $h$, then $(1+h_n/n)^n$ converges to $e^h$.
The above estimates imply
$$\lim_{n\to\infty}-n(1-\hat D(k_n))\,\big(A(k_n)+B(k_n)\big)^{-1}\,{B(k_n)}(1-\hat D(k_n))^{-1}=-K_\alpha\,|k|^\am$$
and
$$\lim_{n\to\infty}\frac{A(0)}{A(k_n)+B(k_n)}=1.$$
We thus have proved Theorem \ref{thm-endpoint} subject to Lemma \ref{lemmaError} and Proposition \ref{propMain2}.
We want to emphasize that the bounds on the lace expansion coefficients $\pi_n(x)$ enter the calculation only through (\ref{eqMain2}) and the error bound in Lemma \ref{lemmaError}.

\subsection{Bounding the lace expansion coefficients}\label{sect_LaceExpansionCoefficients}\noindent
In this section we prove an estimate on moments of the lace expansion coefficients $\pi_n(x)$. This estimate is used to prove Proposition \ref{propMain2}.
Let us begin by stating the moment estimate.
\bl[Finite moments of the lace expansion coefficients]\label{lemmaMoments}
For $\alpha>0$, $d>2(\am)$ and $L$ sufficiently large,
we let
\begin{equation}\label{eqDefDelta}
    \delta
    \begin{cases}
        \in\big(0\,,\,(\am)\wedge(d-2(\am))\big)\qquad&\text{if $\alpha\neq2$,}\\
        =0\qquad&\text{if $\alpha=2$.}
    \end{cases}
\end{equation}
Then, for any $z\le z_c$,
\begin{equation}\label{eqLemmaMoments}
    \sum_{x\in\Zd}\sum_{n=0}^\infty|x|^{(\am)+\delta}\,|\pi_n(x)|\,z^n<\infty.
\end{equation}
\el
The fact that the $((\am)+\delta)$th moment of $\Pi_{z_c}(x)$ exists is the key to the proof of \eqref{eqMain2}.
Interestingly, there is a crossover between the phases $\alpha<2$ and $\alpha>2$, with $\alpha=2$ playing a special role.
A version of Lemma \ref{lemmaMoments} in the setting of oriented percolation is contained in \cite[Proposition 3.1]{ChenSakai09}.

Before we start with the proof of Lemma \ref{lemmaMoments}, we shall review some basic facts about structure and convergence of quantities related to $\pi_n(x)$ introduced in \eqref{eqExpansionX}--\eqref{eqDefPi}.
Our main reference for that is the monograph by Slade \cite{Slade06}, who gives a detailed account of the lace expansion for self-avoiding walk.
Other references are \cite{Hofst05,MadraSlade93}.
We shall also need results from \cite{HeydeHofstSakai08}, where a long-range version of the step distribution is considered.
For $n\ge2$, $N\ge1$, $x\in\Zd$, there exist quantities $\pi_n^\ssup{N}(x)\ge0$ such that
\begin{equation}\label{eqPiBd1}
    \pi_n(x)=\sum_{N=1}^\infty (-1)^N\pi_n^\ssup{N}(x).
\end{equation}
A combination of Theorem 4.1 with Lemma 5.10 (both references to Slade \cite{Slade06}), together with $\beta=O(L^{-d})$ \cite[Prop.\ 2.2]{HeydeHofstSakai08} shows
\begin{equation}\label{eqPiBd2}
    \sum_{x\in\Zd}\sum_{n=2}^\infty  \pi_n^\ssup{N}(x)\,z_c^n
    <O(L^{-d})^N,
\end{equation}
where the constant in the $O$-term is uniform for all $N$.
Consequently, \eqref{eqPiBd2} is summable in $N\ge1$ provided that $L$ is sufficiently large, and hence
\begin{equation}\label{eqPiBd3}
    \hat\Pi_{z_c}(k)
    \le \sum_{x\in\Zd}\sum_{n=2}^\infty |\pi_n(x)|\,z_c^n
    <\infty.
\end{equation}

Lemma \ref{lemmaMoments} implies Proposition \ref{propMain2}, as we will show now.
\proof[Proof of Proposition \ref{propMain2} subject to Lemma \ref{lemmaMoments}.]
We first prove the assertion for $\alpha\le2$, and afterwards consider $\alpha>2$.

For $\alpha\le2$, we choose $\delta\ge0$ satisfying (\ref{eqDefDelta}) and such that $\alpha+\delta\le2$.
Then we use $0\le1-\cos(k\cdot x)\le |k\cdot x|^{\alpha+\delta}$ to estimate
\begin{eqnarray}
  \nonumber
  \left|\hat\Pi_{z_c}(0)-\hat\Pi_{z_c}(k)\right|
  &\le& \sum_{x\in\Zd}\sum_{n=2}^\infty[1-\cos(k\cdot x)]\,|\pi_n(x)|\,z_c^n \\
  \nonumber
  &\le&  \sum_{x\in\Zd}\sum_{n=2}^\infty
       |k\cdot x|^{\alpha+\delta}\,|\pi_n(x)|\,z_c^n \\ \label{eqMomentCalc1}
  &\le&  |k|^{\alpha} \,|k|^{\delta} \sum_{x\in\Zd}\sum_{n=2}^\infty |x|^{\alpha+\delta}\,|\pi_n(x)|\,z_c^n.
\end{eqnarray}
We use \eqref{eqPropD5} and Lemma \ref{lemmaMoments} to bound further
\begin{equation}\label{eqMomentCalc2}
    \frac{|\hat\Pi_{z_c}(0)-\hat\Pi_{z_c}(k)|}{1-\Dk}=
        \begin{cases}
            O(|k|^\delta)\quad&\text{if $\alpha<2$},\\
            O(1/\log(1/|k|))\quad&\text{if $\alpha=2$,}
        \end{cases}
\end{equation}
which proves (\ref{eqMain2}) for $\alpha\le2$.

For $\alpha>2$, we fix $\delta\in(0,2\wedge(d-4))$.
We apply the Taylor expansion
\begin{equation}\label{eqMomentCalc3}
    1-\cos(k\cdot x)=\frac12(k\cdot x)^2+O(|k\cdot x|^{2+\delta}),
\end{equation}
together with spatial symmetry of the model and Lemma \ref{lemmaMoments} to obtain
\begin{equation}\label{eqMomentCalc4}
    \hat\Pi_{z_c}(0)-\hat\Pi_{z_c}(k)
    = \sum_{x\in\Zd}\sum_{n=2}^\infty[1-\cos(k\cdot x)]\,\pi_n(x)\,z_c^n
    = \frac{|k|^2}{2d}\sum_{x\in\Zd}\sum_{n=2}^\infty|x|^2\,\pi_n(x)\,z_c^n
      +O(|k|^{2+\delta}).
\end{equation}
Eq.\ \eqref{eqMain2} for $\alpha>2$ now follows from \eqref{eqMomentCalc4} and \eqref{eqPropD5}.
\qed
\vspace{.5em}

In the remainder of the section we prove Lemma \ref{lemmaMoments}.
A key point in the proof is the use of a new form of (spatial) fractional derivative, first applied by Chen and Sakai \cite{ChenSakai09} in the context of oriented percolation.

\proof[Proof of Lemma \ref{lemmaMoments}.]
For $t>0$, $\zeta\in(0,2)$,
we let
\begin{equation}\label{eqDefK}
    K'_\zeta:=\int_0^\infty\frac{1-\cos(v)}{v^{1+\zeta}}\,\d v
    \in(0,\infty),
\end{equation}
yielding
\begin{equation}\label{eqMom1}
    t^\zeta=\frac1{K'_\zeta}\int_0^\infty\frac{1-\cos(ut)}{u^{1+\zeta}}\,\d u.
\end{equation}

For $\alpha>0$ and $d>2(\am)$, we choose $\delta$ as in (\ref{eqDefDelta}).
For $x\in\Zd$ we write $x=(x_1,\dots,x_d)$. Then by reflection and rotation symmetry of $\pi_n(x)$,
\begin{equation}\label{eqMom2}
    \sum_{x\in\Zd}\sum_{n=0}^\infty|x|^{(\am)+\delta}\,|\pi_n(x)|\,z^n
    \le d^{\,{((\am)+\delta)}/2+1}\sum_{x\in\Zd}\sum_{n=0}^\infty |x_1|^{(\am)+\delta}\,\sum_{N=2}^\infty\pi_n^{\sss (N)}(x)\,z_c^n,
\end{equation}
cf.\ \cite[Lemma 4.1]{ChenSakai09}.
We now apply (\ref{eqMom1}) with $\zeta=\delta_1,\delta_2$, given by
\begin{eqnarray}
  \delta_1 &\in&
    \big(\delta\,,\,(\am)\wedge(d-2(\am))\big),\\
  \delta_2 &=& (\am)+\delta -\delta_1.
\end{eqnarray}
This yields
\begin{equation}\label{eqMom4}
    O(1)\int_0^\infty\frac{\d u}{u^{1+\delta_1}}\int_0^\infty\frac{\d v}{v^{1+\delta_2}}
    \sum_{x\in\Zd}\sum_{n=0}^\infty \sum_{N=2}^\infty [1-\cos(u\,x_1)]\,[1-\cos(v\,x_1)]\,\pi_n^{\sss (N)}(x)\,z_c^n
\end{equation}
as an upper bound of (\ref{eqMom2}).
We write the double integral appearing in \eqref{eqMom4} as the sum of four terms, $I_1+I_2+I_3+I_4$, where
\begin{equation}\label{eqMom3}
    I_1=
    \sum_{N=2}^\infty\int_0^1\frac{\d u}{u^{1+\delta_1}}\int_0^1\frac{\d v}{v^{1+\delta_2}}
    \sum_{x\in\Zd}\sum_{n=0}^\infty  [1-\cos(\uv\cdot x)]\,[1-\cos(\vv\cdot x)]\,\pi_n^{\sss (N)}(x)\,z_c^n
\end{equation}
with
\begin{equation}\label{eqDefUvVv}
\uv=(u,0,\dots,0)\in\R^d,\qquad\vv=(v,0,\dots,0)\in\R^d,
\end{equation}
and $I_2$, $I_3$, $I_4$ are defined similarly:
\begin{equation}\label{eqMom5}
    I_2=\int_0^1\d u\int_1^\infty\d v\cdots,\quad
    I_3=\int_1^\infty\d u\int_0^1\d v\cdots,\quad
    I_4=\int_1^\infty\d u\int_1^\infty\d v\cdots.
\end{equation}
We now show that $I_1,\dots,I_4$ are all finite, which implies \eqref{eqLemmaMoments}.
The bound $I_4<\infty$ simply follows from $1-\cos t\le2$ and \eqref{eqPiBd3}.
In order to prove the bounds $I_1,I_2,I_3<\infty$ we need the particular structure of the $\pi_n^{\sss (N)}(x)$-terms.

To this end, we define
\begin{equation}\label{eqDefGTsaw}
    \tilde G_z(x)=z(D\ast G_z)(x),
    \qquad x\in\Zd,
\end{equation}
and
\begin{equation}\label{eqDefBT}
    \tilde B(z)=\sup_{x\in\Zd}(G_z\ast \tilde G_z)(x).
\end{equation}
In \cite[Theorem 4.1]{Slade06} it is shown that for $z\ge0$, $N\ge1$,
\begin{equation}\label{eqDiaBoundSlade1}
    \sum_{x\in\Zd}[1-\cos(k\cdot x)]\,\Pi_z^\ssup{1}(x)=0
\end{equation}
and
\begin{equation}\label{eqDiaBoundSlade2}
    \sum_{x\in\Zd}[1-\cos(k\cdot x)]\,\Pi_z^\ssup{N}(x)
    \le \frac N2(N+1)\left(\sup_x\,[1-\cos(k\cdot x)]\,G_z(x)\right) \tilde B(z)^{N-1},
    \quad N\ge2.
\end{equation}
These bounds are called \emph{diagrammatic estimates}, because the lace expansion coefficients $\pi_z^\ssup{N}(x)$ are expressed in terms of diagrams, whose structure is heavily used in the derivation of the above bounds.
The composition of the diagrams and their decomposition into two-point functions as in \eqref{eqDiaBoundSlade1}--\eqref{eqDiaBoundSlade2} is described in detail in \cite[Sections 3 and 4]{Slade06}.
It is clear that a slight modification of this procedure proves the bound
\begin{equation}\label{eqDiagrammaticBoundSAW}
\begin{split}
    & \sum_{x\in\Zd}\sum_{n=0}^\infty  [1-\cos(\vv\cdot x)]\,[1-\cos(\uv\cdot x)]\,\pi_n^{\sss (N)}(x)\,z^n \\
    & \quad\le\, O(N^4)\,
    \tilde B(z)^{N-2}\,
    \bigg(\sup_x\,[1-\cos(\vv\cdot x)]\,G_z(x)\bigg)\\
    &\quad\qquad\times
    \bigg(\sup_y\sum_{x\in\Zd}\,[1-\cos(\uv\cdot x)]\,G_z(x)\,G_z(y-x)\bigg).
\end{split}
\end{equation}
Given \eqref{eqDiagrammaticBoundSAW}, it remains to show the following three bounds:
\begin{eqnarray}
   \tilde B(z_c)=\sup_{x\in\Zd}(G_{z_c}\ast \tilde G_{z_c})(x) &\le& O\big(L^{-d}\big); \label{eqDiaBound1}\\
   \sup_x\,[1-\cos(\vv\cdot x)]\,G_{z_c}(x) &\le&
   O\!\left(v^\am\right); \label{eqDiaBound2}\\
   \sup_y\sum_{x\in\Zd}\,[1-\cos(\uv\cdot x)]\,G_{z_c}(x)\,G_{z_c}(y-x)&\le& O\!\left(u^{\left(d-2(\am)\right)\wedge(\am)}\right).\label{eqDiaBound3}
\end{eqnarray}
Suppose \eqref{eqDiaBound1}--\eqref{eqDiaBound3} were true, then
\begin{equation}
\begin{split}
    &\sum_{x\in\Zd}\sum_{n=0}^\infty  [1-\cos(\uv\cdot x)]\,[1-\cos(\vv\cdot x)]\,\pi_n^{\sss (N)}(x)\,{z_c}^n\\
    &\qquad\le O\!\left(N^4\right)O\big(L^{-d}\big)^{N-2}\,O\!\left(v^{\am}\right)O\!\left(u^{\left(d-2(\am)\right)\wedge(\am)}\right).
\end{split}
\end{equation}
Since $\delta_1<(\am)\wedge(d-2(\am))$ and $\delta_2<\am$, we obtain that $I_1$ is finite for $L$ sufficiently large, as desired.
Similarly, it follows that $I_2$ and $I_3$ are finite.
It remains to prove \eqref{eqDiaBound1}--\eqref{eqDiaBound3}, and we use results from \cite{HeydeHofstSakai08} to prove it.

We introduce the quantity
\begin{equation}\label{eqDefLambdaZ}
\lambda_z:=1-\frac{1}{\GN}=1-\frac{1}{\chi(z)}\in[0,1].
\end{equation}
Then $\lambda_z$ satisfies the equality
\begin{equation}
    \GN=\ClN,
\end{equation}
where $\ClK=[1-\lambda_z\Dk]^{-1}$ is the Fourier transform of the simple random walk Green's function.
This definition is motivated by the intuition that $\Gk$ and $\ClK$ are comparable in size and, moreover, the discretized second derivative
\begin{equation}
    \Delta_k\hat G_z(l):=\hat{G}_z(l-k)+\hat{G}_z(l+k)-2\hat G(l)
\end{equation}
is bounded by
\begin{equation}\label{eqDefU}
    \Ul:=200\,\ClK^{-1}\left\{ \hat{C}_{\lambda_z}(l-k)\hat{C}_{\lambda_z}(l)+\hat{C}_{\lambda_z}(l)\hat{C}_{\lambda_z}(l+k)+\hat{C}_{\lambda_z}(l-k)\hat{C}_{\lambda_z}(l+k)\right\}.
\end{equation}

To make this more precise, we consider the function $f\colon[0,z_c]\to\R$, defined by
\begin{equation}\label{eqDefF}
f:=f_1\vee f_2\vee f_3
\end{equation}
with
\begin{equation}\label{eqDefF1}
f_1(z):=z,
\qquad f_2(z):=\sup_{k\in\Td}\frac{\Gk}{\ClK},
\end{equation}
and
\begin{equation}\label{eqDefF3}
f_3(z):=\sup_{k,l\in\Td}\frac{|\Delta_k\hat G_z(l)|}{\Ul},
\end{equation}
It is an important result in \cite{HeydeHofstSakai08} that, under the conditions of Theorem \ref{thm-endpoint}, the function $f$ is uniformly bounded on $[0,z_c)$, cf.\ \cite[Prop.\ 2.5 and 2.6]{HeydeHofstSakai08}.
In fact, it is shown that $f(z)\le1+O(L^{-d})$, but for our need it suffices to have $f$ uniformly bounded.
Since the bound is uniform, 
we can conclude that even $f(z_c)<\infty$.

Indeed, \eqref{eqDiaBound1} follows by standard methods from \cite[Proposition 2.2]{HeydeHofstSakai08}, see e.g.\ \cite[(5.28) in conjunction with Lemma 5.10]{Slade06}.
Furthermore, \eqref{eqDiaBound2} is proven in \cite[Lemma B.5]{HeydeHofstSakai08} in the context of the Ising model, but applies verbatim to self-avoiding walk.
It remains to prove \eqref{eqDiaBound3}.
Since
\begin{eqnarray}
    \nonumber&&
    \sup_y\sum_{x\in\Zd}\,[1-\cos(\uv\cdot x)]\,G_{z_c}(x)\,G_{z_c}(y-x)\\ \nonumber
    &=& \sup_y\int_{\Td} \e^{-il\cdot y}\left(\G_{z_c}(l)-\frac12\left(\G_{z_c}(l-\uv)+\G_{z_c}(l+\uv)\right)\right)\G_{z_c}(l)\;
    \frac{\operatorname{d}\!l}{(2\pi)^{d}}\\
    &\le&\int_{\Td} \left|\frac12\,\Delta_\uv\, \G_{z_c}(l)\right|\G_{z_c}(l)\;
    \frac{\operatorname{d}\!l}{(2\pi)^{d}}\;,
\end{eqnarray}
our bounds $f_2(z_c)\le K$ and $f_3(z_c)\le K$, together with ${\lambda_{z_c}}=1$, imply that
\begin{equation}
    \begin{split}
    &\sup_y\sum_{x\in\Zd}\,[1-\cos(\uv\cdot x)]\,G_{z_c}(x)\,G_{z_c}(y-x)\\
    &\quad\le100K^2\, \hat{C}_1(\uv)^{-1} \int_{\Td} \left(\hat{C}_1(l-\uv)\,\hat{C}_1(l+\uv)+\hat{C}_1(l-\uv)\,\hat{C}_1(l)\right.\\
    &\hspace{5cm}\left.{}+\hat{C}_1(l)\,\hat{C}_1(l+\uv)\right) \hat{C}_1(l) \; \frac{\operatorname{d}\!l}{(2\pi)^{d}}\\
    &\quad=O(1)\,[1-\hat D(\uv)]\int_{\Td} \left(
        \frac{1}{[1-\hat D(l-\uv)]\,[1-\hat D(l+\uv)]\,[1-\hat D(l)]}\right.\\
        &\hspace{1.5cm}\left.{}+
        \frac{1}{[1-\hat D(l-\uv)]\,[1-\hat D(l)]^2}+
        \frac{1}{[1-\hat D(l+\uv)]\,[1-\hat D(l)]^2}
    \right)\frac{\operatorname{d}\!l}{(2\pi)^{d}}. \label{eqCosGxBound1}
    \end{split}
\end{equation}
Chen and Sakai show that the integral term on the right hand side of \eqref{eqCosGxBound1} is bounded above by $O\big(u^{(d-3(\am))\wedge0}\big)$, cf.\ \cite[(4.30)]{ChenSakai09}.
Furthermore, $1-\hat D(\uv)\le O\big(u^\am\big)$ by \eqref{eqPropD5}.
The combination of the above inequalities implies \eqref{eqDiaBound3}, and hence the claim follows.
\qed
\vspace{.5em}

\subsection{Error bounds}\label{sectErrorBound}\noindent
The proof of Lemma \ref{lemmaError} is the final piece in the proof of Theorem \ref{thm-endpoint}.
Our proof of Lemma \ref{lemmaError} makes use of the following lemma:
\bl\label{lem-PowerTails}
Consider a function $g$ given by the power series  $g(z)=\sum_{n=0}^\infty a_nz^n$, with $z_c$ as radius of convergence.
\begin{enumerate}
  \item If $|g(z)|\le O(|z_c-z|^{-b})$ for some $b\ge 1$, then $|a_n|\le O(z_c^{-n}\,\log(n))$ if $b=1$, or $|a_n|\le O(z_c^{-n}\,n^{b-1})$ if $b>1$.
  \item If $|g'(z)|\le O(|z_c-z|^{-b})$ for some $b>1$, then $|a_n|\le O(z_c^{-n}\,n^{b-2})$.
\end{enumerate}
\el
The proof of assertion $(i)$ is contained in \cite[Lemma 3.2]{DerbeSlade98}, and $(ii)$ is a direct consequence of $(i)$ since $(i)$ implies that $|n\,a_n|\le O(z_c^{-n}\,n^{b-1})$.
Lemma \ref{lem-PowerTails} is the key to the proof of Lemma \ref{lemmaError}.
\proof[Proof of Lemma \ref{lemmaError}]
We recall
\begin{equation}\label{eqProof10a}
    \Theta_z(k)=\sum_{n=0}^\infty\theta_n(k)\,z^n,
\end{equation}
where
\begin{equation}\label{eqProof8a}
    \Theta_z(k)=\frac{[1-z/z_c]\,E_z(k)}
                    {\big([1-z/z_c]\left(A(k)+E_z(k)\right)+B(k)\big)\,
                     \big([1-z/z_c]\,A(k)+B(k)\big)}.
\end{equation}
We fix $\eps\in(0,(d(\am)^{-1}-2)\wedge1)$ and aim to prove $|\theta_n(k)|\le O(z_c^{-n}\,n^{-\eps})$, where the constant in the $O$-term is uniform for $k\in\Td$.
By Lemma \ref{lem-PowerTails} it is sufficient to show $|\partial_z\Theta_z(k)|\le O\big(|z_c-z|^{-(2-\eps)}\big)$.

Before bounding $\partial_z\Theta_z(k)$, we consider derivatives of $\hat\Pi_z(k)$ (the Fourier transform of $\Pi_z(x)$ introduced in \eqref{eqDefPi}).
The first derivative of $\partial_z\hat\Pi_{z}(k)$ is  converging absolutely for $z\le z_c$, i.e.,
\begin{equation}\label{eqPiBd4}
    \sum_{x\in\Zd}\sum_{n=2}^\infty n\,|\pi_n(x)|\,z_c^{n-1}
    <\infty,
\end{equation}
cf.\ \cite[Theorem 6.2.9]{MadraSlade93} for a proof in the finite-range setting, and again \cite{HeydeHofstSakai08} for the extension to long-range systems.
Moreover, we claim that
\begin{equation}\label{eqPiBd5}
    \sum_{x\in\Zd}\sum_{n=2}^\infty n(n-1)^\eps\,|\pi_n(x)|\,z_c^{n-1}
    <\infty;
\end{equation}
for $\eps\in(0,(d(\am)^{-1}-2)\wedge1)$.
The bound \eqref{eqPiBd5} can be proved by considering \emph{temporal} fractional derivatives, as introduced in \cite[Section 6.3]{MadraSlade93}.
In particular, the proof of \cite[Theorem 6.4.2]{MadraSlade93}
shows
\begin{equation}\label{eqPiBd5a}
    \sup_{x\in\Zd}\sum_{n=2}^\infty n(n-1)^\eps\,c_n(x)\,z_c^{n-1}
    \le O(1)\int_{\Td}\sum_{n\ge2}n(n-1)^\eps\Dk^{n-2}\dk,
\end{equation}
(see the first displayed identity in
\cite[p.\ 196]{MadraSlade93}).
On the one hand, \eqref{eqPropD3} and \eqref{eqPropD5} imply that there exists some constant $c_1>0$ such that $1-\Dk\ge c_1\,|k|^\am$ for all $k\in\Td$, whence
$    \Dk=1-(1-\Dk)
    \le \e^{-(1-\Dk)}
    \le \e^{-c_1\,|k|^\am}$.
On the other hand, $-\Dk\le1-c_2$ for a positive constant $c_2$, by \eqref{eqPropD4}.
Together these bounds yield
\begin{eqnarray}
    \int\limits_{\Td}\Dk^{n-2}\dk
    &\le& \int\limits_{\substack{k\in\Td\colon\\ \Dk\ge0}}
            \e^{-c_1\,(n-2)\,|k|^\am}\dk\nnb
    &&  {}+ \int\limits_{\substack{k\in\Td\colon\\ \Dk<0}}
            (1-c_2)^{n-2}\dk\nnb
    &\le& O(n^{-d/(\am)})+(1-c_2)^{n-2}\le O(n^{-d/(\am)}).
\end{eqnarray}
Hence the right hand side of \eqref{eqPiBd5a} is less than or equal to
\begin{equation}\label{eqPiBd5b}
    \sum_{n\ge2}n(n-1)^\eps\,O(n^{-d/(\am)}),
\end{equation}
and this is finite if $1+\eps-d/(\am)<-1$.
Furthermore, the proof of \cite[Corollary 6.4.3]{MadraSlade93} shows that
\begin{equation}
    \sum_{x\in\Zd}\sum_{n=2}^\infty n(n-1)^\eps\,|\pi_n(x)|\,z_c^{n-1}
    \le O(1)\,\left(\sup_{x\in\Zd}\sum_{n=2}^\infty n(n-1)^\eps\,c_n(x)\,z_c^{n-1}\right)
\end{equation}
under the conditions of Theorem \ref{thm-endpoint}.
This proves \eqref{eqPiBd5}.

We now prove that
\begin{equation}\label{eqErrorProof5}
    E_z(k)\le O(|z_c-z|^\eps)
\end{equation}
by considering the power series representation of $\hat\Pi_z(k)$ in \eqref{eqProof5}:
\begin{equation}\label{eqErrorProof6}
    E_z(k)=\frac1{z_c-z}\sum_x\sum_{n\ge2}\e^\ikx\pi_n(x)\left(z_c^n-z^n\right)
            -\sum_x\sum_{n\ge2}\e^\ikx\pi_n(x)\,n\,z_c^{n-1}.
\end{equation}
Since
\begin{equation}\label{eqErrorProof7}
    \frac{z_c^n-z^n}{z_c-z}=\sum_{l=0}^{n-1} z^l\;z_c^{(n-1)-l},
\end{equation}
one has
\begin{equation}\label{eqErrorProof8}
    E_z(k)=\sum_x\sum_{n\ge2}\e^\ikx\pi_n(x)\sum_{l=1}^{n-1}\left(z^l-z_c^l\right)z_c^{(n-1)-l}.
\end{equation}
For every $\zeta,\eps\in(0,1)$ and $n\ge2$,
\begin{eqnarray}\nonumber
    \left|1-\zeta^{n-1}\right|
    &=&   \left|(1-\zeta^{n-1})^{1-\eps}\,\left(\frac{1-\zeta^{n-1}}{1-\zeta}\right)^{\eps}\,(1-\zeta)^{\eps}\right|\\
    \label{eqErrorProof11}
    &\le &\left|\sum_{l=0}^{n-2}\zeta^l\right|^\eps\,(1-\zeta)^{\eps}
    \le (n-1)^\eps\left(1-\zeta\right)^\eps.
\end{eqnarray}
Applying this for $\zeta=z/z_c$, we obtain for $z<z_c$ and $0<l<n$,
\begin{eqnarray}\nonumber
    \left|z^l-z_c^l\right|z_c^{(n-1)-l}
    &=& \left|1-\left(\frac{z}{z_c}\right)^l\right|z_c^{n-1}
    \le \left|1-\left(\frac{z}{z_c}\right)^{n-1}\right|z_c^{n-1}\\
    \label{eqErrorProof12}
    &\le& \left|1-\frac{z}{z_c}\right|^\eps (n-1)^\eps\,z_c^{n-1}.
\end{eqnarray}
Insertion into \eqref{eqErrorProof8} yields
\begin{equation}\label{eqErrorProof13}
    \left|E_z(k)\right|
    \le\left(z_c-z\right)^\eps\,\sum_x\sum_{n\ge2}n(n-1)^\eps\,|\pi_n(x)|\,z_c^{n-1}
    \le O(|z_c-z|^{\eps}),
\end{equation}
where the last bound uses \eqref{eqPiBd5}.
We further differentiate \eqref{eqProof5} to get
\begin{eqnarray}\nonumber
    \partial_z E_z(k)
    &=& \frac{(z_c-z)\;\partial_z\big(\hat\Pi_{z_c}(k)-\hat\Pi_{z}(k)\big)+\big(\hat\Pi_{z_c}(k)-\hat\Pi_{z}(k)\big)}{(z_c-z)^2}\\
    &=& \frac1{z_c-z}\left(\frac{\hat\Pi_{z_c}(k)-\hat\Pi_{z}(k)}{z_c-z}-\partial_z\hat\Pi_z(k)\right).
\end{eqnarray}
A calculation similar to \eqref{eqErrorProof6}--\eqref{eqErrorProof13} shows
\begin{equation}\label{eqErrorProof10}
    \left|\partial_z E_z(k)\right|
    \le \left|\frac{E_z(k)}{z_c-z}\right|
        +\frac{1}{z_c-z}
        \left|\sum_x\sum_{n\ge2}\e^\ikx\pi_n(x)\,n\left(z_c^{n-1}-z^{n-1}\right)\right|
    \le O(|z_c-z|^{\eps-1}).
\end{equation}

We write $D_1$ and $D_2$ for the two factors in the denominator in \eqref{eqProof8a}. Then
\begin{equation}\label{eqErrorProof3}
\begin{split}
    z_c^2\,\partial_z\Theta_z(k)
    =&\frac {z_c}{D_1\,D_2} \,
        \big((z_c-z)\,\partial_zE_z(k)-E_z(k)\big) \\
        &{}- \frac {z_c-z}{(D_1\,D_2)^2}\,E_z(k)\Big(
            \big(-A(k)-E_z(k)+(z_c-z)\,\partial_zE_z(k)\big)\,D_2
            -D_1\,A(k)
        \Big).
\end{split}
\end{equation}
The $D_1$- and $D_2$-term in the numerator in the second line of \eqref{eqErrorProof3} can be canceled with the denominator, so that $D_1$ and $D_2$ appear only in the denominator. It is therefore sufficient to give lower bounds on them.
Indeed, there is a constant $c>0$ such that
\begin{equation}\label{eqErrorProof14}
    |D_1|
    =\left|z_c\,\Gk\right|^{-1}
    \ge z_c^{-1}\,\chi(z)
    \ge c\left(z_c-z\right),
\end{equation}
where the last bound follows from \cite[(1.24) and Theorem 1.3]{HeydeHofstSakai08}.
Furthermore, there are constants $c',C>0$ such that
\begin{equation}\label{eqErrorProof15}
    |D_2|
    \ge |D_1|-|E_z(k)\left(z_c-z\right)|
    \ge c\left(z_c-z\right)-C\left(z_c-z\right)^{1+\eps}
    \ge c'\left(z_c-z\right),
\end{equation}
by \eqref{eqErrorProof8} and \eqref{eqErrorProof14}.
The lower bounds on $D_1$ and $D_2$, together with the bounds on $E_z(k)$ and $\partial_zE_z(k)$ in \eqref{eqErrorProof5} and \eqref{eqErrorProof10}, prove that
\eqref{eqErrorProof3} is uniformly bounded for all $z\le z_c$, and in particular
\begin{equation}\label{eqErrorProof4}
    |{\partial_z}\Theta_z(k)|
    \le O(|z_c-z|^{-(2-\eps)}).
\end{equation}
Finally, assertion $(ii)$ in Lemma \ref{lem-PowerTails} implies
\begin{equation}
    |\theta_n(k)|\le O(z_c^{-n}\, n^{-\eps})
\end{equation}
for all $\eps\in(0,(d(\am)^{-1}-2)\wedge1)$, uniformly in $k$.
\qed
\vspace{.5em}

\section{The mean-$r$ displacement: Proof of Theorem \ref{thm-xiR}}\label{6s:RadiusOfGyration}\noindent
\proof[Proof of Theorem \ref{thm-xiR}.]
We start the proof by noting that the reflection and rotation symmetry of $c_n$ implies
\begin{equation}\label{eqMoment1}
    \frac{1}{c_n}\sum_{x\in\Zd}|x|^rc_n(x)
    \asymp \sum_{x\in\Zd}|x_1|^r\;\frac{c_n(x)}{c_n},
\end{equation}
where $x_1$ denotes the first component of the vector $x\in\Zd$.
Recalling \eqref{eqDefFAlpha}, it is therefore sufficient to prove
\begin{equation}\label{eqMoment2}
    \sum_{x\in\Zd}|x_1|^r\;\frac{c_n(x)}{c_n}
    \asymp f_\alpha(n)^{-r}.
\end{equation}
The upper and lower bound in \eqref{eqMoment2} are proved separately, by different methods. We start with the former.

Our proof of the upper bound uses methods similar to those developed in Section \ref{sect_LaceExpansionCoefficients}, and again a key ingredient is the equality in \eqref{eqMom1}.
Again, we denote by $\uv$ the vector $\uv=(u,0,\dots,0)\in\Rd$.
We consider the generating function of the left hand side of \eqref{eqMoment2},
\begin{equation}\label{eqMoment5}
    H_{z,r}:=\sum_{x\in\Zd}\sum_{n=0}^\infty|x_1|^r c_n(x)z^n,
\end{equation}
and claim that $H_{z,r}\le O(1)\,(z_c-z)^{-1-r/(\am)}$ for $\alpha\neq2$ and
$H_{z,r}\le O(1)\,(z_c-z)^{-1-r/2}\,\log(z_c-z)^{-1/2}$ for $\alpha=2$.
Indeed, by \eqref{eqMom1},
\begin{eqnarray}
    H_{z,r}
    &=&
    \sum_{x\in\Zd}\sum_{n=0}^\infty\int_0^\infty
        \frac{\d u}{u^{1+r}}\,
        [1-\cos(\uv\cdot x)]\,c_n(x)\,z^n\nnb
    &\le&
    \int\limits_0^{(z_c-z)^{1/(\am)}}
        \hspace{-1em}\frac{\d u}{u^{1+r}}\;
        \Big(\hat G_z(0)-\hat G_z(\uv)\Big)
    +
    \int\limits_{(z_c-z)^{1/(\am)}}^\infty
        \frac{\d u}{u^{1+r}}\;
        2\,\GN
        \label{eqMoment6}
\end{eqnarray}
where in the last integral we bounded $1-\cos t\le2$.
The generating function $\Gk$ near the critical threshold $z_c$ is known to be bounded by $O(z_c-z)^{-1}$, cf.\ \cite[Theorem 1.1]{HeydeHofstSakai08} (the ansatz in \eqref{eqProof7} leads to the same bound).
Hence the second integral in \eqref{eqMoment6} is bounded above by
\begin{equation}\label{eqMoment7}
        \int\limits_{(z_c-z)^{1/(\am)}}^\infty
        2\,\GN\frac{\d u}{u^{1+r}}
    \le \frac{O(1)}{z_c-z}        \int\limits_{(z_c-z)^{1/(\am)}}^\infty
        \frac{\d u}{u^{1+r}}
    =   \frac{O(1)}{(z_c-z)^{1+r/(\am)}}.
\end{equation}
The first integral on the right of \eqref{eqMoment6} can be expressed as
\begin{equation}\label{eqMoment8}
    \int\limits_0^{(z_c-z)^{1/(\am)}}
        \hspace{-1em}\frac{\d u}{u^{1+r}}\;
        \GN\,\hat G(\uv)\,\Big(z\big(1-\hat D(\uv)\big)+\big(\hat\Pi_z(0)-\hat\Pi_z(\uv)\big)\Big).
\end{equation}
The proof of Proposition \ref{propMain2} might be extended straightforwardly to show
$$\hat\Pi_z(0)-\hat\Pi_z(\uv)=C_\alpha\big(1-\hat D(\uv)\big)+o(1)\,\big(1-\hat D(\uv)\big)$$ for a certain constant $C_\alpha\ge0$ (with $C_\alpha=0$ if $\alpha\le2$), and the $o(1)$-term vanishes as $u\to0$.
Consequently, \eqref{eqMoment8} is bounded above by
\begin{equation}\label{eqMoment9}
    \frac{O(1)}{(z_c-z)^2}
        \int\limits_0^{(z_c-z)^{1/(\am)}}
        \hspace{-1em}\frac{1-\hat D(\uv)}{u^{1+r}}\,\d u.
\end{equation}

Suppose for now that $\alpha\neq2$,
then $1-\hat D(\uv)\le O(u^\am)$ by \eqref{eqPropD5},
and \eqref{eqMoment9} becomes
\begin{equation}\label{eqMoment10}
        \frac{O(1)}{(z_c-z)^2}
        \int\limits_0^{(z_c-z)^{1/(\am)}}
        \hspace{-.6em}{u^{(\am)-(1+r)}}\;\d u.
    =       \frac{O(1)}{(z_c-z)^{1+r/(\am)}}.
\end{equation}
Consequently, $H_{z,r}\le(z_c-z)^{-1-r/(\am)}$, and Lemma \ref{lem-PowerTails}(i) may be applied to deduce
$$  \sum_{x\in\Zd}|x_1|^r\,c_n(x) \le n^{\frac{r}{\am}}\,z_c^{-n}.$$
An application of Corollary \ref{corol-CnAsymptotic} then finishes the proof of the upper bound in \eqref{eqMoment2}.

If on the other hand $\alpha=2$, then \eqref{eqPropD5} and \eqref{eqMoment9} obtain
\begin{equation}\label{eqMoment11}
    H_{z,r}\le\frac{O(1)}{(z_c-z)^{1+r/2}}\,\log(z_c-z)^{-1/2}.
\end{equation}
We then apply the following version of Lemma \ref{lem-PowerTails}(i) (which may be proved along the same lines as Lemma \ref{lem-PowerTails}):
If $\sum_na(n)\,z^n\le(z_c-z)^{-b}\log(z_c-z)^{-1/2}$ for some $b>1$, then $|a(n)|\le O(1)\, n^b\,\log n^{1/2}$.
Together with Corollary \ref{corol-CnAsymptotic} this obtains
$$ \sum_{x\in\Zd}|x_1|^r\,\frac{c_n(x)}{c_n}\le n^{r/2}\;\log\sqrt{n}\qquad\text{for $\alpha=2$}.$$

Finally, we complement the proof of the theorem by showing the lower bound in \eqref{eqMoment2}.
It follows from Theorem \ref{thm-endpoint} that
\begin{equation}\label{eqMoment3}
    \lim_{n\to\infty}1-\frac{\hat c_n(\uv_n)}{\hat c_n(0)}
    =1-\exp\{-K_\alpha\,|u|^\am\},
\end{equation}
and the limit is strictly positive as long as $u\neq0$.
Hence there exists a positive constant $b=b(d,\alpha,L)$ such that for $u=1$ and all $n\in\N$,
\begin{equation}\label{eqMoment4}
    b
    \:\le\:1-\frac{\hat c_n(\uv_n)}{\hat c_n(0)}
    \:=\:\sum_{x\in\Zd}\big[1-\cos\big(u\,f_\alpha(n)\,x_1\big)\big]\,\frac{c_n(x)}{c_n}
    \:\le\:\sum_{x\in\Zd}f_\alpha(n)^r\,|x_1|^r\,\frac{c_n(x)}{c_n},
\end{equation}
where we used $1-\cos t\le|t|^r$ for $r\le\am$ in the last bound.
This implies the lower bound in \eqref{eqMoment2}, and proves the theorem.
\qed

\section{Convergence of finite dimensional distributions:\\ Proof of Theorem \ref{prop-FinDimConv}}\label{6s:FinDimConvergence}\noindent

\proof[Proof of Theorem \ref{prop-FinDimConv}]
The proof is via induction over $N$, and is very much inspired by the proof of \cite[Theorem 6.6.2]{MadraSlade93}, where finite-range models were considered.
The flexibility in the last argument of $n\bT$ is needed to perform the induction step.
We shall further write $nt^{\sss (j)}$ and $nT$ instead of $\lfloor nt^{\sss (j)}\rfloor$ and $\lfloor nT\rfloor$ for brevity.

To initialize the induction we consider the case $N=1$.
Since $\cnT{1}_{n\bT}(\bk_n)=\hat c_{nT}(k^{\sss(1)}_n)$, the assertion for $N=1$ is a minor generalization of Theorem \ref{thm-endpoint}.
In fact, if we replace $n$ by $nT$, then instead of \eqref{eqMain1} we have
\begin{equation}\label{eqFinDimProof0}
    nT\,[1-\hat D(k_n)]
    =nt^{\sss (1)}(1-g_n)\left[1-\hat D\big(f_\alpha(t^{\sss (1)}n)\,k\,(t^{\sss (1)})^{1/(\am)}\big)\right]
    \rightarrow
    |k|^\am\,t^{\sss (1)}
    \quad\text{as $n\to\infty$}.
\end{equation}
With an appropriate change in \eqref{eqProof12} we obtain \eqref{eqFinDimConv} for $N=1$ from Theorem \ref{thm-endpoint}.

To advance the induction we prove \eqref{eqFinDimConv} assuming that it holds when $N$ is replaced by $N-1$.
For an $n$-step walk $w\in\Wcal_n$ and $0\le a\le b\le n$ it will be convenient to write
\begin{equation}\label{eqDefKw}
K_{[a,b]}(w):=\1_{\{(w_a,\dots,w_b)\text{ is self-avoiding}\}}.
\end{equation}
We further consider the quantity $J_{[a,b]}(w)$ that arises in the algebraic derivation of the lace expansion as in \cite[Sect.\ 3.2]{Slade06}. For our needs it suffices to know that
\begin{equation}\label{eqDefJpi}
    \sum_{w\in\Wcal_n(x)}W(w)J_{[0,n]}(w)=\pi_n(x)
\end{equation}
and, for any integers $0\le m\le n$ and $w\in\Wcal_n$,
\begin{equation}\label{eqKJKexpansion}
    K_{[0,n]}(w)=\sum_{I\ni m}K_{[0,I_1]}(w)\,J_{[I_1,I_2]}(w)\,K_{[I_2,n]}(w),
\end{equation}
where the sum is over all intervals $I=[I_1,I_2]$ of integers with either $0\le I_1<m<I_2\le n$ or $I_1=m=I_2$.
We refer to \cite[(3.13)]{Slade06} for \eqref{eqDefJpi},
and to \cite[Lemma 5.2.5]{MadraSlade93} for \eqref{eqKJKexpansion}.
By \eqref{eqDefCnN2} and \eqref{eqKJKexpansion},
\begin{equation}\label{eqFinDimProof1}
    \cnT{N}_{n\bT}(\bk_n)
    =\sum_{I\ni nt^{\sss (N-1)}}
        \sum_{w\in \Wcal_{nT}}
        \e^{i\bk_n\cdot\Delta w(n\bT)}
        W(w)\;
        K_{[0,I_1]}(w)\,J_{[I_1,I_2]}(w)\,K_{[I_2,nT]}(w).
\end{equation}
Let $\cnTl{N}$ and $\cnTg{N}$ denote the contributions towards \eqref{eqFinDimProof1} corresponding to intervals $I$ with length $|I|=I_2-I_1\le b_n$ and $|I|>b_n$, respectively.
It will turn out that the latter contribution is negligible.
We take $n$ sufficiently large so that $(nt^{\sss (N-1)}-nt^{\sss (N-2)})\vee(nt^{\sss (N)}-nt^{\sss (N-1)})\ge b_n$ and
\begin{equation}\label{eqFinDimProof2}
\begin{split}
    &\cnTl{N}_{n\bT}(\bk_n)
    =\sum_{\substack{I\ni \,nt^{\sss (N-1)}\\|I|\le b_n}}
        \cnT{N-1}_{(nt^{\sss(1)},\dots,nt^{\sss(N-2)},I_1)}
        \big(k_n^{\sss (1)},\dots,k_n^{\sss (N-1)}\big)\;
        \times\;\hat c_{nT-I_2}(k_n^{\sss (N)})\\
        &\quad{}\times\!\!\sum_{w\in \Wcal_{|I|}}\!
        \exp\!{\left\{ik_n^{\sss (N-1)}\cdot\, w_{nt^{\sss(N-1)}-I_1}
                +ik_n^{\sss (N)}\cdot\, (w_{I_2-I_1}-w_{nt^{\sss(N-1)}-I_1})\right\}}\;
        W(w)\,
        J_{[0,|I|]}(w).
\end{split}
\end{equation}
We use $\e^y=1+O(|y|^{\alpha\wedge1})$ and \eqref{eqDefJpi} to see that the second line in \eqref{eqFinDimProof2} is equal to
\begin{equation}\label{eqFinDimProof3}
    \sum_x
        \left(1+O(|f_\alpha(n)\,x|^{\alpha\wedge1})\right) \pi_{|I|}(x).
\end{equation}
By the induction hypothesis,
\begin{equation}\label{eqFinDimProof4}
\begin{split}
    &\cnT{N-1}_{(nt^{\sss(1)},\dots,nt^{\sss(N-2)},I_1)}\big(k_n^{\sss (1)},\dots,k_n^{\sss (N-1)}\big)\\
    &\qquad=\hat c_{I_1}(0)\left(
     \exp\Big\{-K_\alpha\,\sum_{j=1}^{N-1}|k^{\sss (j)}|^\am\;(t^{\sss (j)}-t^{\sss (j-1)})\Big\}+o(1)\right)
\end{split}
\end{equation}
and
\begin{equation}\label{eqFinDimProof5}
    \hat c_{nT-I_2}(k_n^{\sss (N)})
    =\hat c_{nT-I_2}(0)\left(
     \exp\Big\{-K_\alpha\,|k^{\sss (N)}|^\am\;(t^{\sss (N)}-t^{\sss (N-1)})\Big\}+o(1)\right),
\end{equation}
where the error terms are uniform in $|I|\le b_n$.

Substituting \eqref{eqFinDimProof3}--\eqref{eqFinDimProof5} into \eqref{eqFinDimProof2} yields
\begin{equation}\label{eqFinDimProof6}
    \cnTl{N}_{n\bT}(\bk_n)
    =\exp\left\{-K_\alpha\,\sum_{j=1}^{N}|k^{\sss (j)}|^\am\;(t^{\sss (j)}-t^{\sss (j-1)})\right\}
    \cnTl{N}_{n\bT}(\bz)+\Theta+o(1)
\end{equation}
where
\begin{equation}\label{eqFinDimProof7}
    |\Theta|\le
    \sum_{\substack{I\ni nt^{\sss (N-1)}\\|I|\le b_n}}
    \hat c_{I_1}(0)\,\hat c_{nT-I_2}(0)\,
    \sum_xO\!\left(|f_\alpha(n)\,x|^{\alpha\wedge1}\right) \pi_{|I|}(x).
\end{equation}
In \eqref{eqFinDimProof7} there are precisely $m-1$ ways to choose the interval $I\ni nt^{\sss (N-1)}$ of length $|I|=m$.
We further bound
\begin{eqnarray}
    \frac{|\Theta|}{\hat c_{nT}(0)}
    &\le&
    \sum_{m=1}^{b_n} m
    \sum_xO\!\left(|f_\alpha(n)\,x|^{\alpha\wedge1}\right) \pi_{m}(x)\,z_c^m
    \nnb
    \label{eqFinDimProof8}
    &\le& O(|f_\alpha(n)|^{\alpha\wedge1}\,b_n)
    \sum_{m=1}^{\infty}
    \sum_x|x|^{\am}\,|\pi_{m}(x)|\,z_c^m
    =o(1),
\end{eqnarray}
where Corollary \ref{corol-CnAsymptotic} is used in the first inequality,
$m\le b_n$ in the second,
and the last estimate uses \eqref{eqDefBn} and Lemma \ref{lemmaMoments}.
Recalling $\cnT{N}_{n\bT}(\bk)=\cnTl{N}_{n\bT}(\bk)+\cnTg{N}_{n\bT}(\bk)$,
\begin{equation}\label{eqFinDimProof9}
    \frac{\cnTl{N}_{n\bT}(\bk_n)}{\hat c_{nT}(0)}
    =\exp\left\{-K_\alpha\,\sum_{j=1}^{N}|k^{\sss (j)}|^\am\;(t^{\sss (j)}-t^{\sss (j-1)})\right\}
    \left(1-\frac{\cnTg{N}_{n\bT}(\bz)}{\hat c_{nT}(0)}\right)
    +\frac{|\Theta|}{\hat c_{nT}(0)}
    +\frac{\cnTg{N}_{n\bT}(\bk_n)}{\hat c_{nT}(0)},
\end{equation}
and it suffices to show ${\cnTg{N}_{n\bT}(\bk_n)}/{\hat c_{nT}(0)}=o(1)$ as $n\to\infty$.
By bounding $|\e^{i\bk_n\cdot\Delta w(n\bT)}|\le1$ in \eqref{eqFinDimProof1}, and using again \eqref{eqDefJpi} and Corollary \ref{corol-CnAsymptotic},
\begin{equation}\label{eqFinDimProof10}
    \frac{\cnTg{N}_{n\bT}(\bk_n)}{\hat c_{nT}(0)}
    \le
    O(1) \sum_{m=b_n+1}^\infty m \sum_x|\pi_m(x)|\,z_c^m,
\end{equation}
which vanishes as $n\to\infty$ by \eqref{eqPiBd4} and the fact that $b_n\to\infty$ as $n\to\infty$.
We have completed the advancement of the induction, and all error terms occurring are uniform in sequences $g=(g_n)$ that satisfy $0\le g_n\le b_n/n$.
This proves \eqref{eqFinDimConv} for all $N\ge1$.
\qed
\vspace{.5cm}

\section{Tightness}\label{6s:Tightness}\noindent
In this section we prove tightness of the sequence $X_n$, the missing piece for the proof of Theorem \ref{thm-pathconvergence}.
Indeed, tightness is implied by Theorem \ref{thm-xiR} and the following tightness criterion.
\bp[{Tightness criterion \cite{Billi68}}]
    The sequence $\{X_n\}$ is tight in $D([0,1],\Rd)$ if the limiting process $X$ has a.s.\ no discontinuity at $t=1$ and there exist constants $C>0$, $r>0$ and $a>1$
    such that for $0\le t_1<t_2<t_3\le1$ and for all $n$,
    \begin{equation}\label{eqTightness1}
        \left\langle|X_n(t_2)-X_n(t_1)|^{r}\,|X_n(t_3)-X_n(t_2)|^{r}\right\rangle_n
        \le C |t_3-t_1|^a.
    \end{equation}
\ep
This proposition is a slight modification of Billingsley \cite[Theorem 15.6]{Billi68},
where (15.21) is replaced by the stronger moment condition on the bottom of page 128 (both references to Billingsley \cite{Billi68}).

\bc[Tightness]\label{cor-tightness}
The sequence $\{X_n\}$ in \eqref{eqDefXn} is tight in $D([0,1],\Rd)$.
\ec
\proof
We first remark that $\alpha$-stable L\`evy motion indeed has a version without jumps at fixed times, and hence no discontinuity at $t=1$ occurs, see e.g. \cite[Theorem 13.1]{Kalle97}.
Fix $r=\frac34\,(\am)$ (in fact, any choice $r\in((\am)/2,\am)$ is possible).
Again we write $nt$ for $\lfloor nt \rfloor$, for brevity.
The left hand side of \eqref{eqTightness1} can be written as
\begin{equation}\label{eqTightness2}
    \frac{f_\alpha(n)^{2r}}{c_n\,(2dK_\alpha)^{2r/(\am)}}\sum_{w\in\Wcal_n}
    |w(nt_2)-w(nt_1)|^r\,|w(nt_3)-w(nt_2)|^r\,W(w)\,K_{[0,n]}(w),
\end{equation}
where $K_{[0,n]}(w)$ was defined in \eqref{eqDefKw}.
Since
\begin{equation}
    K_{[0,n]}(w)\le K_{[0,nt_1]}(w)\,K_{[nt_1,nt_2]}(w)\,K_{[nt_2,nt_3]}(w)\,K_{[nt_3,n]}(w)
\end{equation}
and, by Corollary \ref{corol-CnAsymptotic},
\begin{equation}
    c_n^{-1}\le O(1)\;c_{nt_1}^{-1}\,c_{nt_2-nt_1}^{-1}\,c_{nt_3-nt_2}^{-1}\,c_{n-nt_3}^{-1},
\end{equation}
we can bound \eqref{eqTightness2} from above by
\begin{equation}
\begin{split}
    &\left\langle|X_n(t_2)-X_n(t_1)|^{r}\,|X_n(t_3)-X_n(t_2)|^{r}\right\rangle_n    \\
    &\quad\le \quad
    O(1)\,f_\alpha(n)^{2r}\,\frac1{c_{nt_2-nt_1}}
    \sum_{w\in\Wcal_{nt_2-nt_1}}|w(nt_2-nt_1)|^r\\
    &\hspace{2cm}\times\,\frac1{c_{nt_3-nt_2}}
    \sum_{w\in\Wcal_{nt_3-nt_2}}|w(nt_3-nt_2)|^r\\
    &\quad= \quad O(1)\,f_\alpha(n)^{2r}\,\left(\xi^{(r)}(nt_2-nt_1)\right)^r\,\left(\xi^{(r)}(nt_3-nt_2)\right)^r.
\end{split}
\end{equation}
By Theorem \ref{thm-xiR} and \eqref{eqDefFAlpha},
\begin{equation}
    \left(\xi^{(r)}(nt^\ast-nt_\ast)\right)^r
    \le O(1)\,f_\alpha(n)^{-r}\,(t^\ast-t_\ast)^{r/(\am)}
\end{equation}
for any $0\le t_\ast<t^\ast\le1$,
so that
\begin{equation}
    \left\langle|X_n(t_2)-X_n(t_1)|^{r}\,|X_n(t_3)-X_n(t_2)|^{r}\right\rangle_n
    \le O(1) \,(t_3-t_1)^{2r/(\am)}
    =   O(1) \,(t_3-t_1)^{3/2}.
\end{equation}
This proves tightness of the sequence $\{X_n\}$.
\qed

\proof[Proof of Theorem \ref{thm-pathconvergence}]
The convergence in distribution in Theorem \ref{thm-pathconvergence} is implied by convergence of finite dimensional distributions and tightness of the sequence $X_n$, see e.g.\ \cite[Theorem 15.1]{Billi68}.
Hence, Theorem \ref{prop-FinDimConv} and Corollary \ref{cor-tightness} imply Theorem \ref{thm-pathconvergence}.
\qed

\appendix
\section{Aymptotics of the step distribution}\label{appendixD}
\proof[Proof of \eqref{eqPropD5}]
We consider separately the cases $\alpha>2$ and $\alpha\le2$.

\paragraph{Case $\alpha>2$.}
We expand
$$ \e^{ik\cdot x}=\exp\Big\{i\sum_{j=1}^dk_jx_j\Big\}
    =1+i\sum_{j=1}^dk_jx_j-\frac12\sum_{j,\ell=1}^dk_j\, k_\ell\, x_j\, x_\ell
      +O(|k\cdot x|^{2+\eps})$$
for $0<\eps<(\alpha-2)\wedge1$.
By reflection symmetry,
$$\sum_{x\in\Zd} \sum_{1\le j\le d} k_j\, x_j\,D(x)=0
\quad\text{and}\quad\sum_{x\in\Zd} \sum_{1\le j < n\le d} k_j\, k_\ell\, x_j\, x_\ell \,D(x)=0.$$
Furthermore, as $D$ is symmetric under rotations by ninety degree,
$$\sum_{x\in\Zd}x_1^2\,D(x)=\sum_{x\in\Zd}x_2^2\,D(x)=\cdots=\frac1d \sum_{x\in\Zd}|x|^2\,D(x),$$
so that
\begin{equation}\label{eqApp1}
	\Dk=\sum_{x\in\Zd}\e^{ik\cdot x}\,D(x)
    =1-\frac{|k|^2}{2d}\sum_{x\in\Zd}|x|^2\,D(x)+O\big(|k|^{2+\eps}\big)\sum_{x\in\Zd}|x|^{2+\eps}\,D(x).
\end{equation}
Setting $\sum_{x\in\Zd}|x|^2\,D(x)=2d\,v_\alpha$ proves the claim.

\paragraph{Case $\alpha\le2$.}
The case $\alpha\le 2$ requires a more elaborate calculation.
This part of the proof is adapted from Koralov and Sinai \cite[Lemma 10.18]{KoralSinai07}, who consider the one-dimensional continuous case.
To this end, we write $f=o(g)$ if $f/g$ vanishes as $|k|\to0$.
We can write $D(x)$ as
\begin{equation}
	D(x)=
	c\,\frac{1+g(x)}{|x|^{d+\alpha}},
\end{equation}
where $c$ is a positive constant and $g$ is a bounded function on $\Rd$ obeying $g(x)\to0$ as $|x|\to0$.
By our assumption, $g$ is rotation invariant for $|x|>M$.
We might limit ourselves to the case $|k|\le1/M$ and split the sum defining $\Dk$ as
\begin{equation}\label{eqApp4}
	\Dk = \sum_{|x|\le M}\e^{ik\cdot x}D(x)+
		 \sum_{M<|x|\le 1/|k|}\e^{ik\cdot x}D(x)+
		  \sum_{1/|k|<|x|}\e^{ik\cdot x}D(x). 
\end{equation}
Denote by $S_1$, $S_2$ and $S_3$ the three sums on the right hand side of \eqref{eqApp4}.
A calculation similar to \eqref{eqApp1} shows
\begin{equation}\label{eqApp5}
	S_1= \sum_{|x|\le M}D(x)+ O\big(|k|^2\big)
		=\sum_{|x|\le M}D(x)+
		 \begin{cases} o\big(|k|^\alpha\big)\quad&\text{if $\alpha<2$,}\\
		 	       o\big(|k|^2\log\frac1{|k|}\big)\quad&\text{if $\alpha=2$.}
		 \end{cases}
\end{equation}
For $S_3$ we substitute $x$ by $y/|k|$ yielding
\begin{equation}
	S_3=|k|^{d+\alpha}
		\sum_{\substack{y\in |k|\Zd\\|y|>1}}
		c\,\frac{1+g(y/|k|)}{|y|^{d+\alpha}}\,\e^{ie_k\cdot y},
\end{equation}
where $e_k=k/|k|$ is the unit vector in direction $k$.
By rotation invariance of $g$ and Riemann sum approximation we obtain
\begin{equation}
	S_3=|k|^{\alpha}\left(
		\int_{|y|\ge1}
		c\,\frac{1+g(y/|k|)}{|y|^{d+\alpha}}\,\e^{iy_1}\d y+o(1)\right),
\end{equation}
with $y_1$ being the first coordinate of the vector $y$
and the error term $o(1)$ vanishing as $|k|\to\infty$.
Finally, the dominated convergence (as $|k|\to\infty$) obtains
\begin{equation}
	S_3=|k|^{\alpha}c
		\int_{|y|\ge1}
		\frac{\e^{iy_1}}{|y|^{d+\alpha}}\,\d y+o\big(|k|^\alpha\big).
\end{equation}

Since $D$ is symmetric, the sum defining $S_2$ can be split as
\begin{equation}\label{eqApp2}
	S_2=	\sum_{M<|x|\le 1/|k|}\left(\e^{ik\cdot x}-1-ik\cdot x\right)D(x)+
		\sum_{M<|x|}D(x)-
		\sum_{1/|k|<|x|}D(x).
\end{equation}
Consider first the last sum.
As before, we substitute $x$ by $y/|k|$, use Riemann sum approximation and finally dominated convergence
to obtain
\begin{equation}
	\sum_{1/|k|<|x|}D(x)
	=|k|^{\alpha+d}\sum_{\substack{y\in |k|\Zd\\|y|>1}}
		c\,\frac{1+g(y/|k|)}{|y|^{d+\alpha}}
	=|k|^{\alpha}c
		\int_{|y|\ge1}
		\frac{\e^{iy_1}}{|y|^{d+\alpha}}\,\d y+o\big(|k|^\alpha\big).
\end{equation}
The second sum on the right of \eqref{eqApp2}, together with the complementary sum in \eqref{eqApp5}, obtains the summand 1 on the left of \eqref{eqPropD5}.
It remains to understand the first sum on the right hand side of \eqref{eqApp2}.
We treat this term with the same recipe as above yielding
\begin{equation}\label{eqApp3}
\begin{split}
  & \sum_{M<|x|\le 1/|k|}\left(\e^{ik\cdot x}-1-ik\cdot x\right)D(x) \\
    & {}	= |k|^\alpha c
	  \int_{|k|M\le|y|\le1}
	  \frac{1+g(y/|k|)}{|y|^{d+\alpha}}\left(y_1^2+O\big(|y_1|^{2+\eps}\big)\right)\,\d
	  y+o\big(|k|^\alpha\big).
\end{split}
\end{equation}
For $\alpha<2$ the integral is uniformly bounded in $k$, and hence the dominated convergence theorem can
be used one more time to obtain the desired asymptotics.
However, if $\alpha=2$ then the dominating contribution towards \eqref{eqApp3} is
\begin{equation}
	|k|^2\int_{|k|M\le|y|\le1}
	  \frac{y_1^2}{|y|^{d+\alpha}}\,\d y
	=\frac{|k|^2}{d}\int_{|k|M\le|y|\le1}
	  \frac{1}{|y|^{d}}\,\d y
	= \const|k|^2\left(\log\frac1{|k|}+\log\frac1M\right).
\end{equation}
Summarizing our calculations, we obtain
\begin{equation}
	 \Dk=\sum_{x\in\Zd}D(x)-v_\alpha|k|^\alpha+o\big(|k|^\alpha\big)
	=1-v_\alpha|k|^\alpha+o\big(|k|^\alpha\big)
\end{equation}
for $\alpha<2$,
and
\begin{equation}
	\Dk
	 =1-v_\alpha|k|^2\log\frac1{|k|}+o\Bigg(|k|^2\log\frac1{|k|}\Bigg)
\end{equation}
for $\alpha=2$,
where $v_\alpha$ is composed of the various integrals arising during the proof.
\qed

\vspace{0.5cm}
\noindent {\bf Acknowledgement.}
Research was carried out while the author was affiliated with Technische Universiteit Eindhoven, and supported by the Netherlands Organization for Scientific Research (NWO).
I am indebted to Akira Sakai, Remco van der Hofstad, and Gordon Slade
for kind support during various stages of this project.
I thank Lung-Chi Chen and a referee for many valuable comments on the manuscript,
and the University of Bath for hospitality during my visit in February 2008.

\def\cprime{$'$}


\begin{thebibliography}{10}

\bibitem{Billi68}
P.~Billingsley.
\newblock {\em Convergence of probability measures}.
\newblock John Wiley \& Sons Inc., New York, 1968.

\bibitem{BrydgSpenc85}
D.~C. Brydges and T.~Spencer.
\newblock Self-avoiding walk in {$5$} or more dimensions.
\newblock {\em Comm. Math. Phys.}, 97(1-2):125--148, 1985.

\bibitem{ChenSakai10}
L.-C. Chen and A.~Sakai.
\newblock Asymptotic behavior of the gyration radius for long-range
  self-avoiding walk and long-range oriented percolation.
\newblock In preparation.

\bibitem{ChenSakai08}
L.-C. Chen and A.~Sakai.
\newblock Critical behavior and the limit distribution for long-range oriented
  percolation. {I}.
\newblock {\em Probab. Theory Related Fields}, 142(1-2):151--188, 2008.

\bibitem{ChenSakai09}
L.-C. Chen and A.~Sakai.
\newblock Critical behavior and the limit distribution for long-range oriented
  percolation. {II}: {S}patial correlation.
\newblock {\em Probab. Theory Related Fields}, 145(3-4):435–--458, 2009.

\bibitem{Cheng00}
Y.~Cheng.
\newblock {\em Long Range Self-Avoiding Random Walks above Critical Dimension}.
\newblock PhD thesis, Temple University, August 2000.

\bibitem{DerbeSlade98}
E.~Derbez and G.~Slade.
\newblock The scaling limit of lattice trees in high dimensions.
\newblock {\em Comm. Math. Phys.}, 193(1):69--104, 1998.

\bibitem{HaraSlade92a}
T.~Hara and G.~Slade.
\newblock Self-avoiding walk in five or more dimensions. {I}. {T}he critical
  behaviour.
\newblock {\em Comm. Math. Phys.}, 147(1):101--136, 1992.

\bibitem{HeydeHofstSakai08}
M.~Heydenreich, R.~van~der Hofstad, and A.~Sakai.
\newblock Mean-field behavior for long- and finite range {I}sing model,
  percolation and self-avoiding walk.
\newblock {\em J. Stat. Phys.}, 132(6):1001--1049, 2008.

\bibitem{Hofst05}
R.~van~der Hofstad.
\newblock Spread-out oriented percolation and related models above the upper
  critical dimension: {I}nduction and superprocesses.
\newblock In {\em Ensaios Matem\'aticos [Mathematical Surveys]}, volume~9,
  pages 91--181. Sociedade Brasileira de Matem\'atica, Rio de Janeiro, 2005.

\bibitem{HofstSlade02}
R.~van~der Hofstad and G.~Slade.
\newblock A generalised inductive approach to the lace expansion.
\newblock {\em Probab. Theory Related Fields}, 122(3):389--430, 2002.

\bibitem{Kalle97}
O.~Kallenberg.
\newblock {\em Foundations of modern probability}.
\newblock Probability and its Applications. Springer-Verlag, New York, 1997.

\bibitem{KoralSinai07}
L.~B. Koralov and Ya.~G. Sinai.
\newblock {\em Theory of probability and random processes}.
\newblock Universitext. Springer, Berlin, second edition, 2007.

\bibitem{MadraSlade93}
N.~Madras and G.~Slade.
\newblock {\em The self-avoiding walk}.
\newblock Probability and its Applications. Birkh\"auser Boston Inc., Boston,
  MA, 1993.

\bibitem{SamorTaqqu94}
G.~Samorodnitsky and M.~S. Taqqu.
\newblock {\em Stable non-{G}aussian random processes}.
\newblock Stochastic Modeling. Chapman \& Hall, New York, 1994.

\bibitem{Slade88}
G.~Slade.
\newblock Convergence of self-avoiding random walk to {B}rownian motion in high
  dimensions.
\newblock {\em J. Phys. A}, 21(7):L417--L420, 1988.

\bibitem{Slade89}
G.~Slade.
\newblock The scaling limit of self-avoiding random walk in high dimensions.
\newblock {\em Ann. Probab.}, 17(1):91--107, 1989.

\bibitem{Slade06}
G.~Slade.
\newblock {\em The Lace Expansion and its Applications}, volume 1879 of {\em
  Lecture Notes in Mathematics}.
\newblock Springer-Verlag, Berlin, 2006.

\bibitem{YangKlein88}
W.-S. Yang and D.~Klein.
\newblock A note on the critical dimension for weakly self-avoiding walks.
\newblock {\em Probab. Theory Related Fields}, 79(1):99--114, 1988.

\end{thebibliography}
\end{document}